\documentclass[final,11pt]{elsarticle}
\usepackage{srcltx}
\usepackage{eurosym}
\usepackage{mathtools}
\usepackage{amsmath}
\usepackage{amsfonts}
\usepackage{amssymb}
\usepackage{amsthm}
\usepackage{graphicx}
\usepackage{mathrsfs}
\usepackage{xcolor}
\usepackage{exscale}
\usepackage{latexsym}
\usepackage[pagewise]{lineno}
\usepackage[colorlinks,plainpages=true,pdfpagelabels,hypertexnames=true,colorlinks=true,pdfstartview=FitV,linkcolor=blue,citecolor=red,urlcolor=black]{hyperref}
\PassOptionsToPackage{unicode}{hyperref}
\PassOptionsToPackage{naturalnames}{hyperref}
\usepackage{enumerate}
\usepackage[shortlabels]{enumitem}
\usepackage{bookmark}
\usepackage{wasysym}
\usepackage{esint}
\usepackage[ddmmyyyy]{datetime}
\usepackage[margin=2cm]{geometry}
\makeatletter
\g@addto@macro\normalsize{%
	\setlength\abovedisplayskip{4pt}
	\setlength\belowdisplayskip{4pt}
	\setlength\abovedisplayshortskip{4pt}
	\setlength\belowdisplayshortskip{4pt}
}
\numberwithin{equation}{section}
\everymath{\displaystyle}
\usepackage[capitalize,nameinlink]{cleveref}
\crefname{section}{Section}{Sections}
\crefname{subsection}{Subsection}{Subsections}
\crefname{condition}{Condition}{Conditions}
\crefname{hypothesis}{Hypothesis}{Conditions}
\crefname{assumption}{Assumption}{Assumptions}
\crefname{lemma}{Lemma}{Lemmas}
\crefname{definition}{Definition}{Definitions}

\crefformat{equation}{\textup{#2(#1)#3}}
\crefrangeformat{equation}{\textup{#3(#1)#4--#5(#2)#6}}
\crefmultiformat{equation}{\textup{#2(#1)#3}}{ and \textup{#2(#1)#3}}
{, \textup{#2(#1)#3}}{, and \textup{#2(#1)#3}}
\crefrangemultiformat{equation}{\textup{#3(#1)#4--#5(#2)#6}}%
{ and \textup{#3(#1)#4--#5(#2)#6}}{, \textup{#3(#1)#4--#5(#2)#6}}%
{, and \textup{#3(#1)#4--#5(#2)#6}}

\Crefformat{equation}{#2Equation~\textup{(#1)}#3}
\Crefrangeformat{equation}{Equations~\textup{#3(#1)#4--#5(#2)#6}}
\Crefmultiformat{equation}{Equations~\textup{#2(#1)#3}}{ and \textup{#2(#1)#3}}
{, \textup{#2(#1)#3}}{, and \textup{#2(#1)#3}}
\Crefrangemultiformat{equation}{Equations~\textup{#3(#1)#4--#5(#2)#6}}%
{ and \textup{#3(#1)#4--#5(#2)#6}}{, \textup{#3(#1)#4--#5(#2)#6}}%
{, and \textup{#3(#1)#4--#5(#2)#6}}

\crefdefaultlabelformat{#2\textup{#1}#3}
\numberwithin{equation}{section}

\newtheorem{theorem} {Theorem}[section]
\newtheorem{proposition}[theorem]{Proposition}
\newtheorem{lemma}[theorem]{Lemma}

\newtheorem{example}[theorem]{Example}
\newtheorem{counter example}[theorem]{Counter Example}
\newtheorem{remark}[theorem] {Remark}
\newtheorem{definition}[theorem] {Definition}

\def\N{\mathbb{N}}
\def\CC{{\rm \kern.24em \vrule width.02em height1.4ex depth-.05ex \kern-.26emC}}

\def\TagOnRight

\def\AA{{it I} \hskip-3pt{\tt A}}

\def\QQ{\rlap {\raise 0.4ex \hbox{$\scriptscriptstyle |$}} {\hskip -0.1em Q}}

\makeatletter
\newcommand{\vo}{\vec{o}\@ifnextchar{^}{\,}{}}
\makeatother

\def\YYint#1#2#3{{\setbox0=\hbox{$#1{#2#3}{\iint}$}
		\vcenter{\hbox{$#2#3$}}\kern-.50\wd0}}

\def\XXint#1#2#3{{\setbox0=\hbox{$#1{#2#3}{\int}$}
		\vcenter{\hbox{$#2#3$}}\kern-.50\wd0}}

\makeatletter
\def\namedlabel#1#2{\begingroup
	\def\@currentlabel{#2}%
	\label{#1}\endgroup
}
\makeatother
\makeatletter
\newcommand{\rmh}[1]{\mathpalette{\raisem@th{#1}}}
\newcommand{\raisem@th}[3]{\hspace*{-1pt}\raisebox{#1}{$#2#3$}}
\makeatother

\newcounter{desccount}
\newcommand{\descitem}[2]{\item[#1]\refstepcounter{desccount}\label{#2}}
\newcommand{\descref}[2]{\hyperref[#1]{\textnormal{\textcolor{black}{}\textcolor{blue}{\bf #2}\textcolor{black}{}}}}

\newcommand{\dref}[2]{\hyperref[#1]{\textcolor{black}{(}\textcolor{blue}{\bf #2}\textcolor{black}{)}}}
\newcommand{\be} {\begin{eqnarray}}
\newcommand{\ee} {\end{eqnarray}}
\newcommand{\Bea} {\begin{eqnarray*}}
	\newcommand{\Eea} {\end{eqnarray*}}
\newcommand{\pa} {\partial}
\newcommand{\re}{\mathbb{R}}

\newcommand{\al} {\alpha}
\newcommand{\rr}{\rightarrow}

\newcommand{\dip}{\displaystyle}
\newcommand{\B} {\beta}
\newcommand{\de} {\delta}

\newcommand{\p}  {\prime}
\newcommand{\e}  {\epsilon}
\newcommand{\De} {\Delta}

\newcommand{\la} {\lambda}
\newcommand{\si} {\sigma}

\newcommand{\f}{\infty}
\newcommand{\R}{\mathbb{R}}

\newcommand{\noi} {\noindent}

\DeclareMathOperator{\dv}{div}

\newcommand{\abs}[1]{\left| #1\right|}

\newcounter{whitney}
\refstepcounter{whitney}

\newcounter{ineqcounter}
\refstepcounter{ineqcounter}
\makeatletter
\def\ps@pprintTitle{%
	\let\@oddhead\@empty
	\let\@evenhead\@empty
	\def\@oddfoot{}%
	\let\@evenfoot\@oddfoot}
\makeatother

\newcommand{\refcheckize}[1]{%
	\expandafter\let\csname @@\string#1\endcsname#1%
	\expandafter\DeclareRobustCommand\csname relax\string#1\endcsname[1]{%
		\csname @@\string#1\endcsname{##1}\wrtusdrf{##1}}%
	\expandafter\let\expandafter#1\csname relax\string#1\endcsname
}
\makeatother

\refcheckize{\cref}
\refcheckize{\Cref}


\makeatletter
\newcommand{\mainsectionstyle}{%
	\renewcommand{\@secnumfont}{\bfseries}
	\renewcommand\section{\@startsection{section}{2}%
		\z@{.5\linespacing\@plus.7\linespacing}{-.5em}%
		{\normalfont\bfseries}}%
}
\makeatother
\usepackage{pgf,tikz}
\usetikzlibrary{arrows}
\usetikzlibrary{decorations.pathreplacing}

\usepackage{xpatch}
\xpatchcmd{\MaketitleBox}{\hrule}{}{}{}
\xpatchcmd{\MaketitleBox}{\hrule}{}{}{}

\date{}

\let\oldbibliography\thebibliography
\renewcommand{\thebibliography}[1]{\oldbibliography{#1}
	\setlength{\itemsep}{0pt}}

\allowdisplaybreaks
\begin{document}
	\begin{frontmatter}
		
		\title{Non existence of the BV regularizing effect  for  scalar conservation laws in several space dimension \vspace*{-0.25cm}}	
		\author[myaddress]{Shyam Sundar Ghoshal}\ead{ghoshal@tifrbng.res.in}
		\author[myaddress]{Animesh Jana}\ead{animesh@tifrbng.res.in }
		\address[myaddress]{Centre for Applicable Mathematics, Tata Institute of Fundamental Research, Post Bag No 6503, Sharadanagar, Bangalore - 560065, India. \vspace*{-1cm}}

		\begin{abstract}
		This article deals with the regularity aspects of entropy solutions to scalar conservation laws. We show that for each $C^2$ flux in multi-D, there exists an entropy solution which does not belong to $BV_{loc}(\R^d)$ for all time. For this purpose, we construct a non-$BV_{loc}$ solution in 1-D for a special class of $C^2$ fluxes whose second derivative has a zero. 
		It covers all the $C^2$ functions for which Lax-Ole\u{\i}nik's BV regularizing result is not applicable and provides a classification of one dimensional $C^2$ fluxes based on $L^\f$-$BV_{loc}$ regularizing of entropy solution. In the later part of this article, we extend our result to fractional Sobolev spaces for a class of non-degenerate fluxes.
		
	\end{abstract}

		\end{frontmatter}
	\vspace*{-1.2cm}
\tableofcontents
\section{Introduction}\label{intro}
We consider the following multi-dimensional scalar conservation laws
\begin{align}
\frac{\pa}{\pa t}u+ \sum\limits_{i=1}^d \frac{\partial }{\partial x_i}F_i(u)&=0,\quad x\in\mathbb{R}^d, \ t>0,\label{eqn:conlaw}\\
u(x,0)&=u_0(x),\quad x\in\mathbb{R}^d,\nonumber
\end{align}
where $F=(F_1,\cdots,F_d)\in C^2(\re,\re^d)$. Henceforth we use $F$ as flux function for $d\geq1$ and $f$ as a flux for $d=1$. Existence and uniqueness of solutions to \eqref{eqn:conlaw} is well known due to Kru\v{z}kov \cite{Kruzkov}. Throughout the article, we deal with entropy solution that is a weak solution $u\in C([0,\f);L^{1}_{loc}(\R^d))$ satisfying Kru\v{z}kov entropy criteria \cite{Kruzkov}. 

Our main focus is to study $L^{\f}$-$BV_{loc}$ regularizing effect for entropy solutions to \eqref{eqn:conlaw}. By $L^\f$-$BV_{loc}$ regularizing we mean that the entropy solution gains $BV_{loc}$ regularity for $t>0$ when initial data $u_0\in L^{\f}(\R^d)$. In one dimension, $L^{\f}$-$BV_{loc}$ regularizing of the entropy solution is known for uniformly convex fluxes due to Lax and Ole\u{\i}nik \cite{Lax,Oleinik}.  
This  $L^\f$-$BV_{loc}$ regularizing fails for some special choice of non uniformly convex fluxes in 1-D \cite{AGG_finer,CJ14,Cheng83,Cheverry,DeWe,GGJJ} and for $C^\f$ fluxes in multi-D \cite{CJ17}. In one dimension, near a vanishing point of $f^{\p\p}$, these fluxes locally behave like ${u}^p$ with $u>0$. For some interval $I\subset\R$, we consider a class of $C^2(I,\R)$ functions such that $\{f^{\p\p}=0\}$ is non-empty. For each flux in this class, we prove the existence of an entropy solution not belonging to $BV_{loc}(\R)$. This class contains fluxes considered in \cite{AGG_finer,CJ14,Cheng83,Cheverry,DeWe,GGJJ}. It is the optimal set of $C^2$ fluxes for obtaining a non-$BV$ propagation because the complement of this class in $C^2$ is the set of uniformly convex or concave functions. For higher dimension i.e. for $d>1$ we extend this result to whole $C^2(I,\R^d)$ that is, for each $C^2$ function $F$ we construct an entropy solution which is not in $BV_{loc}(\R^d)$. Note that previously it was known \cite{CJ17} only for $C^\f$ fluxes whereas our result is valid for each $C^2$ flux. Most of the previously constructed examples are non-$BV_{loc}$ up to a finite time. In this article, we prove that the constructed examples are non-$BV_{loc}$ for all time $t>0$.

Since $L^\f$-$BV_{loc}$ regularizing is not possible one can ask a natural question whether the $L^\f$ entropy solution lies in a space larger than $BV_{loc}$. Lions, Perthame and Tadmor \cite{Lions1} first obtained regularity of entropy solution in fractional Sobolev space $W^{s,1}$ (see section \ref{sec:function_space} for the definition) when flux satisfies the following non-degeneracy condition,
\begin{eqnarray}\label{condition:LPT}
\begin{array}{llll}
\mbox{there exists }\al\in (0,1]\mbox{ and } \ C\geq 0 \mbox{ such that }\\
\mathcal{L}^1\left(\{|v|<R_0,\ |\tau +F^\p(v)\cdot \xi|<\delta\}\right)<C \delta^\al,\\
\mbox{for all}\ \delta\in (0,1), \ (\tau,\xi)\in \mathbb{R}\times \mathbb{R}^d \ \mbox{with}\ \tau^2+|\xi|^2=1.
\end{array}
\end{eqnarray}
Here $\mathcal{L}^1$ denotes the one dimensional Lebesgue measure. Authors showed that if $F$ satisfies condition \eqref{condition:LPT} with $\al\in(0,1]$ then the entropy solution belongs to $W^{s,1}_{loc}$ with $s\in\left(0,{\al}/({\al+2})\right)$ by using averaging lemma \cite{Diperna}. For sufficiently smooth flux, it can be shown \cite{Ja10,Ju14,Silvestre} that (\ref{condition:LPT}) means the following: $\{F^{\p\p}(u),\cdots, F^{(m+1)}(u)\}$ spans $\re^d$ for every $u\in I $ where $m$ is the largest integer less than or equal to $1/\al$. For a  survey and comparison between various flux conditions, we refer to \cite{Ju14}. Regularizing in fractional Sobolev space has been studied by several authors \cite{COW,Gess,JP02}. By introducing a new velocity-averaging lemma, Tadmor and Tao \cite{Tao} proved the regularity of entropy solutions in $W^{s,1}$ for $s<{\al}/({2\al+1})$ with an extra condition on flux. Jabin \cite{Ja10} showed that $u\in W^{s,1}(\mathbb{R}^d)$ for all $s<\al$ when flux satisfying \eqref{condition:LPT} and entropy solution satisfies an Ole\u{\i}nik type inequality
\begin{eqnarray}\label{Jabincondition}
\|t\nabla_x\cdot (F^\p(u(\cdot,t)))\|_{\mathcal{M}^1_{loc}(\mathbb{R}^d)}\leq C(\|u_0\|_{L^\f}+\|u_0\|_{L^1}).
\end{eqnarray}
In this context, we would like to mention that for uniformly convex flux entropy solutions of \eqref{eqn:conlaw} with $d=1$ additionally satisfy the following inequality
\begin{equation}\label{Oleinik_ineq}
u(x+h,t)-u(x)\leq\frac{Ch}{t}\mbox{ for all }h,t>0\mbox{ and some constant }C>0.
\end{equation}
The inequality \eqref{Oleinik_ineq} is known as Ole\u{\i}nik's inequality in literature. If an $L^\f$ entropy solution satisfies \eqref{Oleinik_ineq} then it belongs to $BV_{loc}(\R)$ for all $t>0$ \cite{L06}. Hoff introduced a criteria $\dv F^{\p}(u)\leq 1/t$ to characterize multi-D fluxes in  the context of uniqueness of solutions  \cite{Hoff}. In one dimension, it has been proved \cite{Cheng86, Marconi} that $f^{\p}(u)\in BV_{loc}(\R)$ indeed \eqref{Jabincondition} with $d=1$ for the fluxes with polynomial degeneracy. {Note that inequality \eqref{Jabincondition} can be seen as generalized version of Ole\u{\i}nik inequality, Cheng \cite{Cheng86} and Hoff \cite{Hoff} in multi-dimension.} Combining \cite{Cheng86,Marconi} with Jabin's result \cite{Ja10} it can be checked that in 1-D the  entropy solution belongs to $W^{s,1}$ for any $s<\al$ where $\al$ is obtained from flux condition (\ref{condition:LPT}). This result is optimal due to Theorem \ref{theorem2} of present article (see section \ref{main-results}). We show that for $L^\f$ data and $C^2$ flux with non-degeneracy condition \eqref{condition:LPT}, one cannot expect better than $W^{\al,1}$ regularity of the entropy solution to (\ref{eqn:conlaw}). The optimal regularity of the entropy solution in multi-D is still an open question. We refer \cite{AGG_structure,BGJ,CJ-BVphi,CJJ,Golse} for more results on regularizing effects of entropy solution to \eqref{eqn:conlaw}. For finer properties of $L^\f$ entropy solutions in one dimension we refer an interested reader to \cite{BM16, BM17,CJR}.
Rest of the article is organized as follows. Section \ref{main-results} presents main results of this article. We compare our results with known literature in subsection \ref{sec:past_examples}. We recall some definitions and some known facts from theory of conservation laws in section \ref{sec:function_space} and \ref{sec:kruzkov} respectively. We first construct non-BV examples in section \ref{sec:TV}. In section \ref{sec:Wsp} we prove some facts on flux non-degeneracy and construct non-$W^{s,p}$ examples.
\section{Main results}\label{main-results}
Primary goal of this article is to show the non-existence of $BV$ regularizing effect in multi-dimension. We wish to find an initial data defined on $\R^d$ with $d>1$ such that the associated entropy solution to \eqref{eqn:conlaw} does not belong to $BV_{loc}$ for all $t>0$. When $d=1$, we construct non-$BV_{loc}$ solutions for a class of $C^2$ fluxes. More precisely, we prove the following, 

\begin{theorem}\label{theorem1}
	Let $d>1$ and $F\in C^2(\re,\re^d)$. Then there exists an entropy solution $u(\cdot,t)\in  L^\f(\re^d)$ of (\ref{eqn:conlaw})  such that $u(\cdot,t)\notin BV_{loc}(\re^d)$ for all time $t>0$.
\end{theorem}
%
\begin{theorem}\label{theorem:1d}
	Let $I=[a,b]\subset\R$ for some $a<b$. Let $f\in C^2(I,\re)$ be a one dimensional flux function. Suppose there is a point $a_0\in I$ such that $f^{\p\p}(a_0)=0$. Then there exists an initial data $u_0$ such that the corresponding entropy solution $u(\cdot,t)\in  L^\f(\re)$ of \eqref{eqn:conlaw} does not belong to $BV_{loc}(\R)$ for all time $t>0$.
\end{theorem}
Main difference between Theorem \ref{theorem1} and \ref{theorem:1d} is the following: for $d>1$ case we can always find a {unit vector $\xi$ on sphere $\mathcal{S}^{d-1}\subset\R^d$} to confirm $F^{\p\p}(a_0)\cdot\xi=0$. Therefore, Theorem \ref{theorem1} is true for each $F\in C^2(\R,\R^d)$. But for $d=1$, we need to assume that there is a point where the second derivative of $f$ vanishes. We observe that for a flux $f\in C^2(I,\R)$ with $f^{\p\p}\neq0$ on $I$, Theorem \ref{theorem:1d} is not applicable. From works of Lax \cite{Lax} and Ole\u{\i}nik \cite{Oleinik} we know that the $L^{\f}$-$BV_{loc}$ regularizing is true when $f^{\p\p}\neq0$ on $I$ and $u_0(x)\in I$ for a.e. $x\in\R$.  Therefore, we get a complete classification of $C^2$ fluxes based on $L^\f$-$BV_{loc}$ regularizing of entropy solution to \eqref{eqn:conlaw}.

Our next aim is to investigate the regularizing in fractional Sobolev space $W^{s,p}$ when flux is satisfying a non-degeneracy condition \eqref{condition:LPT}. It should be noted that for a $C^2$ flux $F$ the condition \eqref{condition:LPT} hardly tells about the optimality of $\al$. \textit{For an interval $I\subset \re$ we define $\al_{sup}(F,I)$ as follows:
	\begin{equation}\label{alpha_sup}
	\al_{sup}(F,I):=\sup\{\al\in(0,1]; F \mbox{ satisfies \eqref{condition:LPT} on $I$ for }\al\mbox{ and some }C>0 \}.
	\end{equation}}
Similar notion of $\al_{sup}$ has been previously introduced in \cite[section 1.1]{Ju14}. For a $C^\f$ flux $F$, it can be shown \cite{Ju14} that $\al_{sup}(F,I)=(d_F[I])^{-1}$ where $d_F$ is defined as follows
\begin{equation}\label{condition:Ju14}
d_F[I]:=\sup\limits_{u\in I}\inf\{k\geq1;\,\mbox{span}\{F^{\p\p}(u),\cdots,F^{(k+1)}(u)\}=\R^d\}.
\end{equation}


 For non-degenerate fluxes satisfying (\ref{alpha_sup}) with $\al_{sup}\in(0,1)$ we want to show that if it has an $L^{\f}$-$W^{s,p}_{loc}$ regularizing effect then fractional exponent $s$ is bounded above by $\al_{sup}(F,I)$. For this purpose, we introduce the following `mild regularity assumption'.
 \begin{definition}\label{def:mild_reg}
 Let $F\in C^2(I,\R^d)$ for $d\geq1$. Let $A_{\xi}=\{x\in{I};\,\xi\cdot F^{\p\p}(x)=0\}$ be the zero set of $F^{\p\p}\cdot\xi$ for $\xi\in\mathcal{S}^{d-1}$ and $I\setminus A_{\xi}$ can be written as union of mutually disjoint open intervals $\{I_k^{\xi}\}$, that is,
 \begin{equation}\label{intervals}
 I\setminus A_{\xi}=\bigcup\limits_{k\in \mathscr{I}_{\xi}}I_{k}^{\xi}\mbox{ with }\xi\cdot F^{\p\p}\neq0\mbox{ in }I_{k}^{\xi}
 \end{equation}
 for some index set $\mathscr{I}_{\xi}\subset\N$. We say $F$ satisfies `mild regularity assumption', if the following holds
 \begin{equation}\label{mild_regularity}
 \inf\limits_{\xi\in\mathcal{S}^{d-1}}\inf\{\al_{sup}(\xi\cdot F,I_{k}^{\xi});\,k\in\mathscr{I}_\xi\}=\al_{sup}(F,I).
 \end{equation} 
  \end{definition}
 Observe that for 1-D $C^2$ fluxes with $\#\{f^{\p\p}=0\}<\f$, the condition \eqref{mild_regularity} holds. We note that any $C^\f$ flux with $\al_{sup}(F,I)>0$ satisfies \eqref{mild_regularity} (see Remark \ref{lemma:C-infinity}). In 1-D there are many $C^2$ fluxes with $\#\{f^{\p\p}=0\}=\f$ for which \eqref{mild_regularity} is true (see example \ref{example:infinitely}). The condition \eqref{mild_regularity} guarantees that considering all directional non-degeneracy on intervals $I^\xi_k$ is enough to conclude the non-degeneracy $\al_{sup}(F,I)$. It is an open question whether the condition \eqref{mild_regularity} remains true for all $C^2$ fluxes. We leave this for future works. For fluxes satisfying $\al_{sup}(I,R)>0$ along with \eqref{mild_regularity} we have the following result.
\begin{theorem}\label{theorem2}
	Let $d\geq 1$ and $I=[a,b]\subset\R$ for some $a<b$. Let $F\in C^2(I,\re^d)$ be a flux satisfying the non-degeneracy condition (\ref{alpha_sup}) with $\al_{sup}(F,I)\in(0,1)$. Let $F$ be satisfying \eqref{mild_regularity} as in Definition \ref{def:mild_reg}. Then, there exists an entropy solution $u(\cdot,t)\in  L^\f(\re^d)$ of (\ref{eqn:conlaw})  such that $u(\cdot,t)\notin W_{loc}^{s,1}(\re^d)$ for all $t>0$ with $s>\al_{sup}(F,I)$.
\end{theorem}
Note that in Theorem \ref{theorem2} we assume that $\al_{sup}(F,I)<1$. In one dimension there exists a flux function $f\in C^2(I)$ for some $I\subset\R$ such that $\al_{sup}(f,I)=1$. We make a remark on $L^{\f}$-$BV_{loc}$ regularizing of entropy solution to \eqref{eqn:conlaw} for such fluxes.
\begin{remark}\label{remark:BV_for_alpha=1}
	Let $I=[a,b]$ for some $a<b$. There exists a flux $f\in C^2(I,\R)$ and an initial data $u_0\in L^{\f}(\R)$ such that $\al_{sup}(f,I)=1$ and if $u$ is the entropy solution to \eqref{eqn:conlaw} with flux $f$ and data $u_0$ then we have $u(\cdot,t)\notin BV_{loc}(\R)$ for all $t>0$.
\end{remark}
Proof of Remark \ref{remark:BV_for_alpha=1} can be found in section \ref{sec:Proof_of_Remark}. 

To prove Theorem \ref{theorem2} we calculate Besov norm of entropy solution and utilize the fact $W^{s,1}\subset B^{s,p,\theta}$ with suitable choice of $p,\theta$ where  $B^{s,p,\theta}$ is Besov space (see section \ref{sec:function_space} for definition). Therefore, the result on failure of regularizing effect is actually true in more general space than $W^{s,p}$. In the Besov space set up the non-existence of regularizing from $L^{\f}$ data with fluxes satisfying \eqref{alpha_sup} and \eqref{mild_regularity}, reads as follows,
\begin{proposition}\label{proposition1}
	Let $d\geq 1$ and $R_0>0$.  Let $F\in C^2([-R_0,R_0],\re^d)$ be a flux satisfying the non-degeneracy condition (\ref{alpha_sup}) with $\al_{sup}(F,[-R_0,R_0])\in(0,1)$ and additionally \eqref{mild_regularity}. Then there exists an entropy solution $u(\cdot,t)\in  L^\f(\re^d)$ of (\ref{eqn:conlaw})  such that $u(\cdot,t)\notin B_{loc}^{s,p,\theta}(\R^d)$ for all $t>0$ with $s>\al_{sup}(F,[-R_0,R_0]),\,p\geq1,\theta\geq1$.	
\end{proposition}
Till now the discussion has been restricted for fluxes having only $C^2$ regularity. For fluxes satisfying non-degeneracy condition \eqref{alpha_sup} from Theorem \ref{theorem2} we get an upper bound for regularizing in fractional Sobolev space, more precisely, there exists a data such that if the associated entropy solution $u\in W^{s,p}$ then $s$ must satisfy $s\leq \al_{sup}$. In \cite{Ju14} it has been shown for $C^\f$ fluxes that $\al_{sup}\leq d^{-1}$. By applying Theorem \ref{theorem2} we can conclude that there is an entropy solution, $u$ to \eqref{eqn:conlaw} for $C^\f$ flux such that $u\notin W^{s,p}$ for $s>d^{-1}$. Therefore, for $C^\f$ flux we recover an explicit upper bound on fractional exponent $s$. In our next proposition, we obtain an explicit bound of $s$ for $C^{2,\B}$ fluxes. By $C^{2,\B}$ function we mean a $C^2$ function whose second derivative belongs to Holder space $C^{\B}$ when $0<\B\leq 1$. 
\begin{proposition}\label{proposition2}
	Let $d>1$ and $F\in C^{2,\B}([a,b],\re^d)$ for some $0<\B\leq1$. Then there exists an entropy solution $u(\cdot,t)\in  L^\f(\re^d)$ of (\ref{eqn:conlaw})  such that $u(\cdot,t)\notin W_{loc}^{s,p}(\R^d)$ for all $t>0$ with $s>1/(1+\B),\,p\geq1$. Moreover, we have $u(\cdot,t)\notin B_{loc}^{s,p,\theta}(\R^d)$ for $s>1/(1+\B),\,p\geq1,\theta\geq1$.
\end{proposition}

\subsection{Comparison with previously constructed examples}\label{sec:past_examples}
 In a short note \cite{Cheng83}, Cheng first showed that $L^{\f}$-$BV$ regularizing fails for fluxes like  ${u^3}/{3},{u^4}/{4}$. With an explicit estimation of Besov semi-norm of entropy solution, De Lellis and Westdickenberg \cite{DeWe} proved that in one dimension there exists an entropy solution which does not belong to $W^{\al,1}$ for flux $f(u)=(p+1)^{-1}|u|^{p+1}$ and $\al>p^{-1}$. Cheverry \cite{Cheverry} constructed an entropy solution, $u$ for flux ${u^3}/{3}$ with piece-wise constant data to show that $u(\cdot,t)\notin BV_{loc}$ for $t\in[0,{1}/{6}]$. In reference \cite{AGG_finer}, authors built an entropy solution $u(\cdot,t)\notin BV_{loc}$ for non-uniform convex flux satisfying a certain condition. They used the method of backward construction from control theory \cite{AGG_exact}.
 
 In one dimension, Castelli and Junca \cite{CJ14} constructed an example, $u$ for $C^{\f}$ flux, $f$ such that $u(\cdot,t)\notin BV^{s+\e}(\R)\cap W^{s+\e,1}(\R)$ for all $\e>0$, $t\in[0,T]$ with some $T>0$ and $s=d_f^{-1}$ where $d_f$  is determined as follows,
  \begin{equation*}
 d_f=\sup_{u\in[-R,R]}\inf\left\{j\geq1;\,f^{(j+1)}(u)\neq 0\right\}<+\f ,\ R>0,
 \end{equation*}
 where $f^{(k)}$ denotes the k-th derivative of $f$. The index $d_f$ determines non-linearity of flux function $f$ and if $f$ is analytic functions $d_f^{-1}$ coincides with $\al_{sup}$ defined as in \eqref{alpha_sup} \cite{Ju14}. One can see that $d_f$ is achieved at some point $\tilde{u}$ and in a neighbourhood of $\tilde{u}$, $f$ behaves like power law type flux $\abs{u}^{1+d_f}$ (for more details see \cite[Theorem 6]{CJ14} and \cite[Example 14]{CJ17}). Then authors built the required solution following the method similar to \cite{Cheng83}. By using the one-dimensional example, they \cite{CJ17} constructed a `mono-phase' entropy solution in multi dimension such that it does not belong to $W^{s,1}_{loc}(\R^d)$ for $s>\al$ where $\al$ measures the non-degeneracy of smooth flux $F$.
 
 On a slightly different note, we mention the counterexamples obtained in \cite{ADGG,SSG15} to show the non existence of $BV$ solution with $BV$ data for discontinuous flux. Unlike the situation for continuous case in the set up of discontinuous flux, the $BV$ regularizing is possible \cite{SSG16} near the interface even for non uniformly convex fluxes.

 From the above discussion, we observe that all the non-BV solutions obtained in the previous works are done either for fluxes of type $\abs{u}^{p+1}$ or for $C^{\f}$ fluxes which behave like ${u}^{p+1}$ with $u>0$. The key difference between these examples and the solution constructed in Theorem \ref{theorem1} and \ref{theorem:1d} is the difficulty to find an initial data which gives a non-BV entropy solution for time $t\in[0,t_0]$. Due to the explicit structure of flux $f(u)=\abs{u}^{1+p}$ one natural choice of data is an arrangement of sequences like $\{n^{-\e-p^{-1}}\}$ with suitable spatial distance \cite{AGG_finer,CJ17,Cheverry,DeWe}. One can also use $x\sin(1/x)$ as initial data to get a non-$BV_{loc}$ entropy solution for power law type flux \cite{CJ14,Cheng83}.
 Now consider the situation when $f$ is a $C^2$ function such that $f^{\p\p}$ vanishes at some point. To construct a non-BV entropy solution for a such flux, we need a sequence like $\{n^{-\e-p^{-1}}\}$. In this case, explicit structure of flux is not known and naturally, no such explicit sequence works. To resolve this matter we extract an implicit sequence $\{a_k\}$ around a vanishing point $a_0$ of $f^{\p\p}$ such that it verifies a certain estimate. Once we get $\{a_k\}$ one may think about constructing a data by arranging $\{a_0,a_1,a_0,a_2,a_0,\cdots\}$ with suitable spatial distance. This works if we have 
 \begin{equation}\label{TV_a}
 \sum\limits_{k\geq1}\abs{a_k-a_0}=\f.
 \end{equation} Unfortunately, this may not hold in general. Since our choice of $\{a_k\}$ is implicit we can not conclude whether \eqref{TV_a} is satisfied or not. We get rid of this hurdle by repeating a similar oscillation sufficiently many times.
 
 Finally we would like to point out the fact that examples given in \cite{CJ14,CJ17,Cheng83,Cheverry,DeWe} are non-BV up to a finite time $t_0$ whereas in this article we construct examples such that they are not in $BV_{loc}(\R)$ for all time $t>0$. Although examples in \cite{AGG_finer,GGJJ} are non $BV_{loc}$ for all time, we can not rely on their construction because they used convexity of $f$ on both side of the point where $f^{\p\p}$ vanishes. In our situation we have information about $f$ only on one side of a zero of $f^{\p\p}$. To give more clarity on this matter, let us consider flux $f(u)=u^3$. Since $f$ is convex for $u>0$ and concave for $u<0$ we can not invoke results of \cite{AGG_finer} whereas we can still find an all time non-$BV_{loc}$ entropy solution for $u^3$ by Theorem \ref{theorem:1d}.
 
\section{Function spaces, some basic results and notations}\label{sec:function_space}
Here we recall some basic definitions and notations which will be used later on. 

	Let $\Omega\subset\mathbb{R}^d$ be an open set and $u\in L^1(\Omega)$. Suppose $\mathcal{B}(\Omega)$ denotes the Borel $\sigma$-algebra. We say $u\in BV(\Omega)$ if for each $i=1,\cdots,d$ there exists a finite signed measure $\mu_i:\mathcal{B}(\Omega)\to \mathbb{R}$ such that
		\begin{equation}\label{def:BVderivative}
		\int\limits_{\Omega}u\frac{\partial \varphi}{\partial x_i}{d}x=-\int\limits_{\Omega}\varphi d\mu_i,
		\end{equation}
		holds for all $\varphi\in C_c^{\f}(\Omega)$. For $u\in BV(\Omega)$, we denote total variation of $u$ by $TV(u,\Omega)$ defined as follows:
		\begin{equation*}
		TV(u,\Omega)=\sup\left\{\sum\limits_{i=1}^{d}\int\limits_{\Omega}\frac{\pa\Phi_i}{\pa x_{i}};\,\Phi\in C_c^{\f}(\Omega,\R^d),\|\Phi\|_{L^{\f}(\Omega,\R^d)}\leq1\right\}.
		\end{equation*}
		When $d=1$, for $\Omega=[a,b]$ with $a<b$, the semi-norm $TV(u,\Omega)$ is equivalent to following point-wise definition (for instance see \cite{Leoni}),
		\begin{equation*}
		TV(u,[a,b])=\sup\left\{\sum\limits_{k=1}^{n}\abs{u(x_k)-u(x_{k+1})};\,a\leq x_1<x_2<\cdots<x_{n+1}\leq b\right\}.
		\end{equation*} 
		For $d\geq1$, $BV$ functions can be characterize in the point-wise variant on lines. More precisely, we can have 
		\begin{theorem}[\cite{Leoni}]\label{theorem:BV}
			Let $Q\subset\re^d$ be open and $d\geq1$. For each $j=1,\cdots,d$ and $h>0$ we define $Q_{h,j}:=\left\{x\in Q;\,x+he_j\in Q\right\}$
			where $\{e_j;1\leq j\leq d\}$ is the standard basis of $\re^d$. Let $u\in BV(Q)$ and derivatives of $u$ are finite signed measures $\mu_j$, $j=1,\cdots,d$  determined by (\ref{def:BVderivative}). Then for each $j=1,\cdots,d$ and $h>0$ we have
			\begin{equation*}
			\int\limits_{Q_{h,j}}|u(y+he_j)-u(y)|\,dy\leq\,h|\mu_j|(Q)\mbox{ and }			\lim\limits_{h\rr0^+}\int\limits_{Q_{h,j}}\frac{|u(y+he_j)-u(y)|}{h}\,dy=|\mu_j|.
			\end{equation*} 
			Conversely, if $u\in L^1(Q)$ is such that 
			\begin{equation*}
			\liminf\limits_{h\rr0^+}\int\limits_{Q_{h,j}}\frac{|u(y+he_j)-u(y)|}{h}\,dy<\f,\mbox{ for each }j=1,\cdots,d.
			\end{equation*}
			Then $u\in BV(Q)$.
		\end{theorem}
	Note that if some function has unbounded point-wise variation on a line then by virtue of Theorem \ref{theorem:BV} it can not belong to $BV(\R)$. This fact will be use in proof of Theorem \ref{theorem:BV}. We say a function $v\in BV_{loc}(\Omega)$ if $v\in BV(\Omega_1)$ for all $\Omega_1\subset\subset\Omega$.
	
	 Let $\Omega$ be a smooth domain in $\R^d$. For a vector $h\in\R$ and $u\in L^p(\Omega)$ with $1<p<\f$ let $\De_h (u,x)$ be the difference operator defined as follows, $\De_h(u,x):=u(x+h)-u(x)$. We also use the following notation 
	 \begin{equation*}
	 \De_h(u,x,\Omega):=\left\{\begin{array}{rl}
	 \De_h(u,x)&\mbox{ if }x,x+h\in\Omega,\\
	 0&\mbox{otherwise.}
	 \end{array}\right.
	 \end{equation*}
	 The modulus of smoothness of order $1$ of the $L^p$ function $u$ can be measured as
	 $$
	 \omega(u,t,\Omega)_p:=\sup\{\|\De_h(u,\cdot,\Omega)\|_{L^p(\Omega)};\abs{h}\leq t\}.
	 $$ 
	 For $s\in(0,1), p\in[1,\f),\theta\in[1,\f)$ the Besov space $B^{s,p,\theta}$  consists of $u\in L^p(\Omega)$ such that
	 $$
	 \abs{u}_{B^{s,p,\theta}(\Omega)}:=\left(\int\limits_{0}^{1}\frac{\omega(u,t,\Omega)^\theta}{t^{1+s\theta}}\,dt\right)^{\frac{1}{\theta}}<\f.
	 $$ 
	 We say $v\in B^{s,p,\theta}_{loc}(\Omega)$ if $v\in B^{s,p,\theta}(\Omega_1)$ for each open set $\Omega_1\subset\subset\Omega$. Above definition of Besov space is adapted from \cite{Devore93}. For various properties of $B^{s,p,\theta}(\R^d)$ we refer \cite{Leoni}.
	 
	 For $s\in(0,1)$ and $p\in[1,\f)$ we define fractional Sobolev space $W^{s,p}(\Omega)$ as follows,
	 $$
	 W^{s,p}(\Omega)=\left\{g\in L^p(\re^d); |u|_{W^{s,p}(\Omega)}^p=\int\limits_{\Omega}\int\limits_{\Omega}\frac{|g(x)-g(y)|^p}{|x-y|^{d+sp}}dxdy<\f\right\}.
	 $$ 
	 The space $W^{s,p}_{loc}(\R^d)$ consists of functions $u$ such that $u\in W^{s,p}(\Omega)$ for each bounded  open set $\Omega\subset\R^d$. When $\theta=p$ we have $W^{s,p}(\Omega)=B^{s,p,p}(\Omega)$ for $ 1\leq p<\f$ \cite{Leoni, Triebel}.

Throughout the paper we use following definitions of convex function, box function and $\mathcal{S}^{d-1}$.
	\begin{itemize}
		
		\item 	Let $f:\re\to\re$ be a real valued function. We say $f$ is ``strictly convex" if the following holds:
		$ 
		f(\la x+(1-\la )y)< \la f(x)+(1-\la)f(y) \ \mbox{ for any }\la\in(0,1)\ \mbox{ and }x,y\in \re.
			$ 
		
		\item 	Suppose $x\in\re$ then by $[x]$ we denote the greatest integer which is less than or equal to $x$. This we refer later as ``box function".
		\item The notation $\mathcal{S}^{d-1}$ stands for unit sphere in $\re^d$, i.e., $ \mathcal{S}^{d-1}=\left\{\xi\in\re^d;|\xi|=1\right\}$.
		
		\item In this article, $\{e_j;\,1\leq j\leq d\}$ denotes the standard basis of $\re^d$.
		
		\item For $a_0\in\re$, by `$a\rr a_0\pm$' we mean limits from right and left respectively.
	\end{itemize}
   Finally, we prove an elementary lemma. This is helpful in the proof of Theorem \ref{theorem2}.
   \begin{lemma}\label{lemma:subsequence}
   	Let $\{\de_n\}_{n\geq1}$ be a real number sequence such that $\de_n\rr0$ as $n\rr\f$. Then there exists an increasing sequence $\{q_k\}_{k\geq1}$ and a subsequence $\{\de_{n_k}\}_{k\geq1}$ such that $\abs{\de_{n_{k}}}=k^{-q_k}$ and $q_{k+1}\geq q_{k}+1 $ for all $k\geq1$.
   \end{lemma}
   \begin{proof}
   	We take $q_1:=1$ and $n_1:=\max\{n;|\de_n|<1\}$. Having defined $q_1,\cdots,q_{k}$ and $n_1,\cdots,n_{k}$ we define $n_{k+1}$ as follows $
   	n_{k+1}:=\max\left\{n;|\de_n|\leq (k+1)^{-q_k-1}\right\}$. Set $q_{k+1}:={-(\log|\de_{n_{k+1}}|)}/{\log(k+1)}$. From  the choice of $n_{k+1}$ we get $|\de_{n_{k+1}}|\leq\frac{1}{(k+1)^{q_k+1}}$. Therefore, we  have $q_{k+1}\geq q_k+1$.	This completes the proof.
   \end{proof}
\section{Some basics from Kru\v{z}kov theory}\label{sec:kruzkov}
In this section, we recall some basic results on scalar conservation laws based on finite speed of propagation for $L^\f$ data and Kru\v{z}kov's uniqueness proof for \eqref{eqn:conlaw}. These will be used heavily in section \ref{sec:TV} and \ref{sec:Wsp}. Let us first recall the result from \cite{Kruzkov} in the following proposition.
\begin{proposition}[Kru\v{z}kov, \cite{Kruzkov}]\label{Prop:Kruzkov}
	Let $d\geq1$. Let $u,v\in L^{\f}(\R^d\times\R_+)\cap C([0,\f);L^{1}_{loc}(\R^d))$ be two entropy solution to \eqref{eqn:conlaw} associated with initial data $u_0$ and $v_0$ respectively. Let $M$ be defined as follows
	\begin{equation*}
	M=\max\left\{\abs{F^{\p}(w)};\,\abs{w}\leq\max\{\|u\|_{ L^{\f}(\R^d\times\R_+)},\|v\|_{ L^{\f}(\R^d\times\R_+)}\}\right\}.
	\end{equation*}
	Then we have for all $R>0$ and $t>0$ 
	\begin{equation*}
	\int\limits_{\abs{x}\leq R}\abs{u(x,t)-v(x,t)}\,dx\ \leq\int\limits_{\abs{x}\leq R+Mt}\abs{u_0(x)-v_0(x)}\,dx.
	\end{equation*}	
\end{proposition}
In next lemma we consider a planar data (i.e. function of one variable and independent of other variables) in multi dimension then we show that the associated entropy solution is also planar. This is classical and elementary and based on a change of variable and right choice of test function \cite{JuAssym}. We give a short proof in the Appendix to make the presentation self contained.
\begin{lemma}\label{lemma:planar_multiD}
	Let $w$ is entropy solution of the following one dimensional conservation law
	\begin{equation}
	\begin{array}{rll}
	\pa_t w+\pa_{x_1} f_1(w)&=0&\mbox{ for }t>0\mbox{ and }x_1\in\re, \label{w1}\\
     w(x_1,0)&=w_0(x_1)&\mbox{ for } x_1\in\re.
	\end{array}
	\end{equation}	Let $v:\R^d\times\R_+\rr\R$ and $v_0:\R^d\rr\R$ be defined as follows
	\begin{equation}
	\begin{array}{rll}
	v(x_1,x_2,\cdots,x_d,t)&=w(x_1,t)&\mbox{ for }(x_1,x_2,\cdots,x_d,t)\in \re^d\times\R_+,\label{lemma:kruzkov_v}\\
	v_0(x_1,x_2,\cdots,x_d)&=w_0(x_1)&\mbox{ for }(x_1,x_2,\cdots,x_d)\in \re^d.
	\end{array}
	\end{equation}
	Then $v$ is the entropy solution to (\ref{eqn:conlaw}) with initial data $v_0$ and flux $F=(f_1,\cdots,f_d)$.
\end{lemma}
Next lemma is also elementary but has important role in proofs of our main results. It establishes the relation between entropy solutions of \eqref{eqn:conlaw} corresponding to two fluxes related by an isomorphism. This is also a corollary to Proposition \ref{Prop:Kruzkov} and we postpone the proof till Appendix.
\begin{lemma}\label{lemma:linear_transform}
	Let $v$ be the entropy solution to (\ref{eqn:conlaw}) with flux $F=(f_1,\cdots,f_d)$ and initial data $v_0$. Let $\bar{v}$ be the entropy solution to (\ref{eqn:conlaw}) with flux $\bar{F}=(\bar{f}_1,\cdots,\bar{f}_d)$ and initial data $\bar{v}_0$. Suppose $L:\re^d\rr\re^d$ is an invertible linear map and $\textbf{c}\in\R^d$ be a vector. Assume that $v_0(x)=\bar{v}_0(L(x)+\textbf{c})$ for all $x\in \re^d$ and $\bar{F}=L(F)$. Then we have $v(x,t)=\bar{v}(L(x)+\textbf{c},t)$ for all $(x,t)\in \re^{d}\times(0,\f)$.
\end{lemma} 
\section{Non-existence of $BV$-solution}\label{sec:TV}

\subsection{One dimensional results}
In this section, we wish to show the existence of non-$BV_{loc}$ entropy solutions in 1-D when second derivative of flux has a zero, say $a_0$. For this purpose we construct an initial data by using a special sequence $\{a_k\}$ converging to the point $a_0$. In general, for a $C^2$-flux the choice of such sequence is tricky. Here is our next lemma precisely stating the properties of a sequence required to construct a non-$BV_{loc}$ entropy solution.
\begin{lemma}\label{lemma:inflection}
	Let $f:[a,b]\rr\re$ be a $C^2$ function for $a<b\in\re$. Let $a_0\in(a,b)$ be a point such that $f^{\p\p}(a_0)=0$ and $f^{\p\p}\neq 0$ in $I=(a_0,a_0+r)\subset[a,b]$ for some $r>0$. Then there are sequences $\{p_k\}_{k\geq1}\subset(1,\f)$ and $\{a_k\}_{k\geq1}\subset(a_0,a_0+r)$ such that
	\begin{equation}\label{eq:sequence0}
	\frac{|f^\p(a_k)-f^\p(a_0)|}{|a_k-a_0|^{p_k}}=\frac{1}{k},\,|a_k-a_0|\leq\frac{1}{k}\mbox{ and }\frac{|f^\p(a_k)-f^\p(a_0)|}{|a_k-a_0|}\leq\frac{1}{k^2}\mbox{ for all }k\geq k_0
	\end{equation}
	for some $k_0\in\N$ with $k_0\geq2$.
\end{lemma}
Before we go to proof of Lemma \ref{lemma:inflection} we put some remarks and observation on the $\{p_k\}$ sequence which is related to the degeneracy of $f^{\p\p}$ near $a_0$. 
\begin{remark}\label{remark:lemma}
	\begin{enumerate}[(i)]
		\item\label{Remark-i} For $f(u)=(p+1)^{-1}\abs{u}^{p+1}$ with $p>1$, \eqref{eq:sequence0} implies that $p_k$ must lie between $\max\{(p-1)/2,p-1\}$ and $p$. For instance we can choose 
		\begin{equation*}
		\begin{array}{rlllll}
		p_k&=(p+1)/2&\mbox{and}&a_k=k^{-\frac{2}{p-1}}&\mbox{for}&1<p\leq3,\\
		p_k&=p-1&\mbox{and}&a_k=k^{-1}&\mbox{for}&p>3.
		\end{array}
		\end{equation*}
		It is possible to give a choice of $a_k$ such that the associated $p_k$ converges to $p$ as $k\rr\f$. For example we can take $p_k=p-k^{-1}$ and $a_k=k^{-k}$ for $k\geq (p-1)/2$.
		\begin{proof}
			Note that $f^{\p}(u)=u\abs{u}^{p-1}$ and $a_0=0$. From \eqref{eq:sequence0} we get
			\begin{equation}\label{eqn:ppk}
			\abs{a_k}^{{p-p_k}}=k^{-1},\,\abs{a_k}\leq k^{-1}\mbox{ and }\abs{a_k}^{p-1}\leq k^{-2}.
			\end{equation}
			We observe $\log\abs{a_k}\leq0$ and subsequently, for $k>1$ we have	$({p-p_k})\log\abs{a_k}=-\log k<0$ which implies ${p-p_k}>0$. Hence we obtain $p_k<p$. Now from \eqref{eqn:ppk} we get $k^{\frac{1}{p-p_k}}\geq k\mbox{ and }k^{\frac{1}{p-p_k}}\geq k^{\frac{2}{p-1}}$. This implies $1\geq p-p_k\mbox{ and }\frac{p-1}{2}\geq p-p_k$.	Hence $p_k\geq \max\{(p+1)/{2},p-1\}$. A simple calculation proves rest of the statements in \ref{Remark-i} of Remark \ref{remark:lemma}.
		\end{proof}
		\item\label{Remark-ii} Let $f\in C^2[0,e^{-1}]$ such that $f^{\p\p}(u)=(\log\abs{u})^{-1}$, then  for $a_k=e^{-k^2}$ we have $p_k>1$ and $p_k\rr1+$ as $k\rr\f$.
		\begin{proof}
			Define $g:[0,e^{-1}]\rr\R$ as $g(x):=\int\limits_{0}^{x}f^{\p\p}(y)\,dy$. Fix a $k\geq1$. We observe that for $x\in [0,a_k]$ we have 
			\begin{equation}\label{inequality:log}
			\abs{\log\abs{x}}^{-1}\geq x^{\frac{1}{k}}.
			\end{equation} 
			To see \eqref{inequality:log} consider $z=1/x$ for $x>0$ and we need to show $(\log z)^{-1}>z^{-k^{-1}}$, equivalently, $z^{\frac{1}{k}}>\log z$ for $z\geq e^{k^2}$. Define $h(z)=z^{\frac{1}{k}}-\log z$. By a simple calculation we have 
			$
			h^{\p}(z)=z^{-1}\left[k^{-1}z^{\frac{1}{k}}-1\right]\geq0$ for $z\geq e^{k^2}$.	Hence $h$ is increasing on $[e^{k^2},\f)$. Therefore, we have $	h(z)\geq h(e^{k^2})=e^{k}-k^{2}\geq0$ for $k\geq1$.	This establishes \eqref{inequality:log}. By using the inequality \eqref{inequality:log} for $x\in[0,a_k]$ we obtain
			\begin{equation}\label{g_log}
			g(a_k)\geq\int\limits_{0}^{a_k}y^{\frac{1}{k}}\,dy\geq \frac{k}{k+1}e^{-k-k^2}.
			\end{equation}
			Note that $\abs{\log\abs{x}}^{-1}\leq \abs{\log\abs{a_k}}^{-1}=k^{-2}$ for $x\in[0,a_k]$. Therefore, we infer 
			\begin{equation}\label{eq:g(a_k)}
			g(a_k)\leq k^{-2}e^{-k^2}.
			\end{equation} Hence we have $a_k^{-1}\abs{g(a_k)}\leq k^{-2}$. The first identity in \eqref{eq:sequence0} enforces  $p_k=-k^{-2}\log\abs{kg(a_k)}$. By \eqref{eq:g(a_k)} we have $k\abs{g(a_k)}\leq k^{-1}e^{-k^2}$. Subsequently, we obtain $\log\abs{kg(a_k)}\leq-k^2-\log k$. This implies $-\log\abs{kg(a_k)}\geq k^2+\log k$. Hence $p_k\geq 1+k^{-2}\log k>1$.  By \eqref{g_log} we have $
			\log(kg(a_k))\geq -k-k^2+2\log k-\log(k+1)$. Therefore, we get
			\begin{equation}\label{eq:pklog}
			p_k=-k^{-2}\log(kg(a_k))\leq k^{-1}+1-2k^{-2}\log k+k^{-2}\log(k+1).
			\end{equation}
			Since $ k^{-1}-2k^{-2}\log k+k^{-2}\log(k+1)\rr0$ as $k\rr\f$ and $p_k>1$ for all $k$, from \eqref{eq:pklog} we conclude that $p_k\rr1$.
		\end{proof}
		\item\label{Remark-iii} For fluxes $f\in C^2[0,1]$ such that $f^{\p}(u)=\mbox{exp}\{-u^{-2}\}$, we have $p_k\rr\f$ as $k\rr\f$ for any choice of $\{a_k\}$ verifying \eqref{eq:sequence0}. One example of $\{a_k\}$ satisfying \eqref{eq:sequence0} is $a_k=k^{-1}$.
		\begin{proof}
			Suppose $\{a_k\}$ satisfies \eqref{eq:sequence0}. Then we have $p_k=\frac{\log(ke^{-1/a_k^{2}})}{\log(a_k)}$. This implies
			\begin{equation}\label{pke1}
			p_k=\frac{1}{\log(a_k)}\left[\log k-a_k^{-2}\right].
			\end{equation} 
			Since $a_k\leq k^{-1}$, we have $\log a_k\leq -\log k$ or equivalently, $\log k/\log a_k\geq -1$. Therefore, \eqref{pke1} becomes $p_k\geq -1+a^{-2}_k(-\log a_k)^{-1}$. Let $b_k=a_k^{-1}$. Hence we get $p_k\geq b_k^2/\log b_k-1$. Since $a_k\rr0$ we have $b_k\rr\f$ and subsequently, $b_k^2/\log b_k\rr\f$. This shows $p_k\rr\f$. Note that if we consider $a_k=k^{-1}$ it follows that $f^{\p}(a_k)a_k^{-1}=ke^{-k^2}\leq k^{-2}$ for $k\geq2$. Therefore, $\{a_k\}$ verifies \eqref{eq:sequence0} with $p_k=-1+k^2/\log k$.
		\end{proof} 
	\end{enumerate}
\end{remark}

From above three examples of fluxes we see three different behaviour of $p_k$ sequence. We also observe in these examples that $p_k$ depends on the choice of $a_k$. To give more insight on this matter let us consider nonlinearity exponent $p$ for convex flux $f$ determined as follows
	\begin{equation}\label{power_law_condition}
	p=\inf\left\{q;\,\frac{\abs{f^{\p}(u)-f^{\p}(v)}}{\abs{u-v}^q}>0\mbox{ for }u,v\in[\underline{u},\overline{u}],u\neq v\right\}.
	\end{equation}
In general the exponent $p$ as in \eqref{power_law_condition} may not be achieved and it can be 1. From \ref{Remark-i} of Remark \ref{remark:lemma} we observe that $p_k$ (as in \eqref{eq:sequence0}) may not converge to the non-degeneracy exponent $p$. For a fixed $f\in C^2(I,\R)$ and $a_0\in I$ there are many choices for pairs $(p_k,a_k)$ and in the proof of Lemma \ref{lemma:inflection} we get an implicit choice of $(p_k,a_k)$.

\begin{proof}[Proof of Lemma \ref{lemma:inflection}]
	We define a function $Q:(a_0,a_0+r)\times[1,\f)\rr[0,\f]$ as follows $Q(a,p):={|f^{\p}(a)-f^{\p}(a_0)|}{|a-a_0|^{-p}}$.	Since $f^{\p\p}(a_0)=0$ we have
	\begin{equation*}
	\lim\limits_{a\rr a_0+}Q(a,1)=\lim\limits_{a\rr a_0+}\frac{|f^{\p}(a)-f^{\p}(a_0)|}{|a-a_0|}=|f^{\p\p}(a_0)|=0.
	\end{equation*}
	Hence for each $k\geq2$ we can choose
	\begin{equation}\label{akin}
	a_k\in\left(a_0,a_0+k^{-1}\right)\mbox{ such that }Q(a_k,1)\leq k^{-2}.
	\end{equation}
	Since $f^{\p\p}\neq0$ in $(a_0,a_0+r)$, we have $C_k:=Q(a_k,1)>0$. Therefore, as $|a_k-a_0|<k^{-1}$, we have
	\begin{equation}\label{kp1}
	Q(a_k,p)=|a_k-a_0|^{1-p}Q(a_k,1)\geq k^{p-1}C_k\mbox{ for any }p>1. 
	\end{equation}
	Note that for fixed $a_k\in(a_0,a_0+r)$, $p\mapsto Q(a_k,p)$ is a continuous function. For a fixed $k\geq2$, by \eqref{kp1} we have $Q(a_k,p)\rr\f$ as $p\rr\f$. From \eqref{akin} we get $Q(a_k,1)<k^{-1}$. Therefore, by Intermediate Value Theorem we can choose $p_k$ such that $Q(a_k,p_k)=k^{-1}$. Hence the proof is completed.
%
\end{proof}
\begin{remark}
	Note that Lemma \ref{lemma:inflection} is also valid if $f^{\p\p}\neq0$ in $(a_0-r,a_0)$ instead of $(a_0,a_0+r)$.  
\end{remark}
Existence of sequence $a_k$ sets the ground to construct an initial data to have a non-BV solution for \eqref{eqn:conlaw} with $d=1$. For fixed $T>0$, we first focus on constructing an entropy solution which is not in BV space up to the time $T$. Later we extend this for all time and in multi dimension. 
\begin{proposition}\label{Proposition:1d_uptoT}
	Let $f\in C^2(\re,\re)$. Suppose there is a point $a_0\in\re$ such that $f^{\p\p}(a_0)=0$. Then there exist constants $\la>0$ and $u^*>0$ such that the following holds: for every $T>0$ there exists an initial data $v_0^T$ verifying the following properties: 
	\begin{enumerate}
		\item $\|v^T_0\|_{L^{\f}(\R)}\leq u^*$ and $supp(v^T_0)\subset[-\la T,\la(T+1)]$.
		\item if $v^T$ is the entropy solution to \eqref{eqn:conlaw} with initial data $v^T_0$, then we have 
		$ 
		TV(v^T(\cdot,t),[-\la T,\la(T+1)])=\f\mbox{ for all }t\in[0,T].
		$ 
	\end{enumerate}
\end{proposition}
Note that in Proposition \ref{Proposition:1d_uptoT}, since the construction depends on the parameter $T$, we use the notation $v_0^T,v^T$ but the constants $\la,u^*$ do not depend on $T$, they depend only on flux $f$. 
\begin{proof}[Proof of Proposition \ref{Proposition:1d_uptoT}] Based on the possible linearity of flux $f$ in some interval, we divide the proof into two cases. First one is the case when there is a non-trivial interval where $f^{\p\p}=0$. Second case is considered when $f$ is non-linear everywhere in the sense that there is no non-trivial interval such that $f^{\p\p}$ vanishes identically in it.
	\begin{description}
		\descitem{Case(1):}{linear} Suppose $f^{\p\p}=0$ on $[\underline{a},\overline{a}]$ for some $\underline{a}<\overline{a}$. Let $\{y_k\}_{k\geq1}$ be an increasing positive real number sequence such that $y_k\rr y_\f $ as $k\rr\f$ for some $y_\f \in\R$. Fix $a,b\in(\underline{a},\overline{a})$ such that $a\neq b$. Consider $w_0^{lin}$ defined as follows
		\begin{equation}\label{def:initial_w_lin}
		w_0^{lin}(x):=\left\{\begin{array}{rcl}
		a&\mbox{ if }&x<y_1,\\
		b&\mbox{ if }&y_{2m-1}<x<y_{2m}\mbox{ for }m\geq1,\\
		a&\mbox{ if }&y_{2m}<x<y_{2m+1}\mbox{ for }m\geq1,\\
		b&\mbox{ if }&x>y_\f .
		\end{array}	\right.
		\end{equation}
		Since $f^{\p\p}=0$ on $[\underline{a},\overline{a}]$ we have $f^{\p}(u)=c_1$ and $f(u)=c_1u+c_2$ when $u\in[\underline{a},\overline{a}]$ for some constants $c_1,c_2$. Hence, for initial data $w_0^{lin}$, the entropy solution to \eqref{eqn:conlaw} coincides with the weak solution of the transport equation $\pa_t w^{lin}+c_1\pa_x w^{lin}=0$. Let $w^{lin}$ be the weak solution to the transport equation $\pa_t w^{lin}+c_1\pa_x w^{lin}=0$ with initial data $w_0^{lin}$. 
		Then $w^{lin}$ can be written as $w^{lin}(x,t)=w_0^{lin}(x-c_1t)$. Hence we have for $t>0$
		\begin{equation}\label{structure:sol_w}
		w^{lin}(x,t):=\left\{\begin{array}{rcl}
		a&\mbox{ if }&x<y_1+c_1t,\\
		b&\mbox{ if }&y_{2m-1}+c_1t<x<y_{2m}+c_1t\mbox{ for }m\geq1,\\
		a&\mbox{ if }&y_{2m}+c_1t<x<y_{2m+1}+c_1t\mbox{ for }m\geq1,\\
		b&\mbox{ if }&x>y_\f +c_1t.
		\end{array}	\right.
		\end{equation} 
		Since $w^{lin}(\cdot,t)$ oscillates infinitely many times with oscillation strength $\abs{a-b}$, it has unbounded variation for each $t\geq0$. Fix a $t\geq0$. We can choose $z_k=c_1t+(y_{k}+y_{k+1})/2$ to obtain
		\begin{align}
		TV(w^{lin}(\cdot,t),[y_1+c_1t,y_\f +c_1t])&\geq \sum\limits_{k\geq1}\abs{w^{lin}(z_k,t)-w^{lin}(z_{k+1},t)}\nonumber\\
		&=\sum\limits_{k\geq1}\abs{a-b}=\f.\label{TVblowup:w}
		\end{align}
		Now we wish to construct an initial data $v_0$ which is supported in a compact set of $\R$ and the associated entropy solution to \eqref{eqn:conlaw}, $v(\cdot,t)$ has unbounded variation for $t\in[0,T]$. To construct $v_0$ we multiply previously constructed data $w_0^{lin}$ by a characteristic function $\chi_I$ for some well-chosen interval $I$. Note that $0$ may lie outside the interval  $[\underline{a},\overline{a}]$ and in that case characteristics can have non-linear propagation and it can cancel the oscillations. To avoid such situation, it is enough to  choose $I$ sufficiently large because the equation \eqref{eqn:conlaw} has `finite speed of propagation' property for $L^\f$ data. Set $M_1$ as follows
		\begin{equation}\label{def:M1}
		M_1:=\max\{\abs{f^{\p}(u)};\,\abs{u}\leq r\}\mbox{ where }r:=\max\{ \abs{a},\abs{b}\}.
		\end{equation}
		We consider $R_1:=3M_1T$. Let $I\subset \R$ an interval defined as $I:=[y_1-R_1,y_\f+R_1]$. We define initial data $v_{0,T}^{lin}$ as follows
		\begin{equation}\label{def:v_0lin}
		v_{0,T}^{lin}:=w_0^{lin}\chi_{I}
		\end{equation}
		where $w_0^{lin}$ is defined as in \eqref{def:initial_w_lin}. Recall that $y_1\geq0$. Hence we have 
		\begin{equation}\label{eq:support_lin}
		supp(v^{lin}_{0,T})\subset[-3M_1T,y_{\f}+3M_1T].
		\end{equation}
		Let $v^T_{lin}$ be the entropy solution to \eqref{eqn:conlaw} with initial data $v_0$. Now by Proposition \ref{Prop:Kruzkov} we have
		\begin{equation}\label{linear_case1}
		\int\limits_{[x_1,x_2]}\abs{w^{lin}(x,t)-v^T_{lin}(x,t)}\,dx\leq \int\limits_{[x_1-M_1T,x_2+M_1T]}\abs{w_0^{lin}(x)-v^{lin}_{0,T}(x)}\,dx
		\end{equation}
		 for $t\in[0,T]$. Take $x_1=y_1-M_1T$ and $x_2=y_\f+M_1T$. Then from our choice of $I$, we have 
		$[x_1-M_1T,x_2+M_1T]=[y_1-2M_1T,y_\f+2M_1T]\subset I$.
		\begin{figure}
			\centering
			\begin{tikzpicture}[scale=.7]

			\draw[thick,->] (-10 ,0) -- (5,0) node[anchor=north west] {$x$};
			\draw[thick,->] (-10 ,0) -- (-10,5.5) node[anchor=north east] {$t$};
			
			\draw[color=black] (-6.2 ,0) -- (-4.2,4) ;
			\draw[color=black] (-5 ,0) -- (-3,4) ;
			\draw[color=black] (-4 ,0) -- (-2,4) ;
			\draw[color=black] (-3.2 ,0) -- (-1.2,4) ;
			\draw[color=black] (-2.6 ,0) -- (-0.6,4) ;
			\draw[color=black] (-2.2 ,0) -- (-0.2,4) ;
			\draw[color=black] (0 ,0) -- (2,4) ;

			\draw[color=blue] plot [smooth] coordinates { (-7.5,0) (-7.2,1.5) (-6,3) (-5,4) (-4.5,5)};
			\draw[color=blue] plot [smooth] coordinates { (1.5,0) (2.5,1.5) (2.5,3) (2.3,4) (2.3,5)};

			\draw[thick][] (-5.6,0.6) node[anchor=north] {$b$};
			
			\draw[thick][] (-4.5,0.5) node[anchor=north] {$a$};
			
			\draw[thick][] (-3.5,0.6) node[anchor=north] {$b$};
			
			\draw[thick][] (-2.8,0.5) node[anchor=north] {$a$};
			\draw[thick][] (-2.3,0.6) node[anchor=north] {$b$};
			
			\draw[thick][] (0.6,0.5) node[anchor=north] {$a$};
			\draw[thick][] (-6.8,0.5) node[anchor=north] {$a$};
			\draw[thick][] (-8.5,0.6) node[anchor=north] {$0$};
			\draw[thick][] (2.5,0.55) node[anchor=north] {$0$};
			
			\draw[thick][] (-0.3,1.2) node[anchor=south west, transform canvas={scale=1.5}]{$\textbf{...} $};
			
			\draw[thick][] (-5.26,-0.22) node[anchor=south west, transform canvas={scale=1.5}]{$\textbf{.} $};					
			\draw[thick][] (-4.39,-0.22) node[anchor=south west, transform canvas={scale=1.5}]{$\textbf{.} $};
			\draw[thick][] (-3.59,-0.22) node[anchor=south west, transform canvas={scale=1.5}]{$\textbf{.} $};
			\draw[thick][] (-2.92,-0.22) node[anchor=south west, transform canvas={scale=1.5}]{$\textbf{.} $};
			\draw[thick][] (-2.39,-0.22) node[anchor=south west, transform canvas={scale=1.5}]{$\textbf{.} $};
			\draw[thick][] (-1.99,-0.22) node[anchor=south west, transform canvas={scale=1.5}]{$\textbf{.} $};
			\draw[thick][] (-1.73,-0.22) node[anchor=south west, transform canvas={scale=1.5}]{$\textbf{.} $};
			\draw[thick][] (-0.26,-0.22) node[anchor=south west, transform canvas={scale=1.5}]{$\textbf{.} $};
			\draw[thick][] (0.74,-0.22) node[anchor=south west, transform canvas={scale=1.5}]{$\textbf{.} $};

			\draw[thick][] (-6.2,-.1) node[anchor=north] {$y_1$};
			\draw[thick][] (-5,-.1) node[anchor=north] {$y_{2}$};
			\draw[thick][] (0.0,-.1) node[anchor=north] {$y_\f$};
			\draw[thick][] (-7.5,-.1) node[anchor=north] {$A$};
			\draw[thick][] (1.5,-.1) node[anchor=north] {$B$};
			
			\draw[dashed] (5,3.5) -- (-10,3.5) node[anchor= east] {$T$};
			\end{tikzpicture}
			\caption{The initial data is constructed by positioning $a$ and $b$ in between points $y_k$ alternatively where $\{y_k\}$ is an increasing sequence converging to $y_\f$ and its support lies in $[A,B]$ for some $A<y_1<y_\f<B$. $a$ and $b$ are connected by linear waves starting from $y_k$ for $k\in\N\cup\{\f\}$. Nonlinear waves are formed at $A$ and $B$. Entropy solution consists of infinitely many oscillations of strength $\abs{a-b}$ up to time $T$.}\label{figure:linear}
		\end{figure}
				 As $v_{0,T}^{lin}(x)=w_0^{lin}(x)$ for all $x\in I$, from \eqref{linear_case1} we have (see Figure \ref{figure:linear} for a clear illustration)
		$
		v^T_{lin}(x,t)=w^{lin}(x,t)\mbox{ for }(x,t)\in[y_1-M_1T,y_\f+M_1T]\times[0,T].
		$ 
		From \eqref{TVblowup:w} we know that $w^{lin}(\cdot,t)$ has unbounded variation in the interval $[y_1+c_1t,y_\f+c_1t]$. Since $[y_1+c_1t,y_\f+c_1t]\subset[y_1-M_1T,y_\f+M_1T]$ for all $t\in[0,T]$, we conclude that 
		$ 
		TV(v^T(\cdot,t),[y_1-M_1T,y_\f+M_1T])=\f \mbox{ for all }t\in[0,T].
		$ 
		Due to $y_1\geq0$, we have $[y_1-M_1T,y_\f+M_1T]\subset[-3M_1T,y_\f+3M_1T]$. Hence we have
		\begin{equation}\label{TVblowup_lin}
		TV(v^T(\cdot,t),[-3M_1T,y_\f+3M_1T])=\f \mbox{ for all }t\in[0,T].
		\end{equation}
		
		\descitem{Case(2):}{non-linear} Now we discuss the case when $\mathcal{Z}:=\{f^{\p\p}=0\}$ has empty interior. Since $f\in C^2(\R)$, the set $\R\setminus \mathcal{Z}$ is open. If $\R\setminus\mathcal{Z}$ is empty then it reduces to the \descref{linear}{Case(1)}. Therefore, we assume $\R\setminus\mathcal{Z}$ is non-empty. Hence, it can be written as union of at most countably many disjoint open intervals. We choose one of such open intervals, say $(a_0,b_0)$. Since intervals are disjoint and $(a_0,b_0)\subset\R\setminus\mathcal{Z}$ we have $f^{\p\p}(a_0)=0$ and $f^{\p\p}(y)\neq 0$ for $-\f<a_0<y<b_0\leq\f$. Without loss of generality we may assume $b_0<\f$ and $f^{\p\p}>0$ in $(a_0,b_0)$. To set the ground for the construction in multi-D here we build an initial data in one dimension with following parameters: $\{\e_k\}_{k\geq1},\{a_k\}_{k\geq1}\cup\{a_0\},\{A_k\}_{k\geq1}$ satisfying
		\begin{equation}\label{condition:ekAk}
		\abs{f^{\p}(a_k)-f^{\p}(a_0)}\leq \e_k\leq A_k\mbox{ and }C:=\sum\limits_{k=1}^{\f}A_k+2\e_1<+\f.
		\end{equation}
		We invoke Lemma \ref{lemma:inflection} to get sequences $\{p_k\}$ and $\{a_k\}$ verifying \eqref{eq:sequence0}. From \eqref{eq:sequence0}, observe that $\e_k=\abs{a_k-a_0}^{p_k}$ and $A_k=\abs{f^{\p}(a_k)-f^{\p}(a_0)}/\abs{a_k-a_0}$ satisfy \eqref{condition:ekAk}. Note that from \eqref{eq:sequence0} we have $\e_k^{-1}\abs{f^{\p}(a_0)-f^{\p}(a_k)}=k^{-1}\leq1$, $A_k\leq k^{-2}$ and $\abs{a_k-a_0}^{-1}\geq k$. By using \eqref{eq:sequence0} we obtain
		\begin{align*}
		A_k=\abs{f^{\p}(a_k)-f^{\p}(a_0)}\abs{a_k-a_0}^{-1}&=(\abs{f^{\p}(a_k)-f^{\p}(a_0)}\abs{a_k-a_0}^{-p_k})\abs{a_k-a_0}^{p_k-1}\\
		&=k^{-1}\abs{a_k-a_0}^{-1}\e_k\geq\e_k.
		\end{align*} Hence $\{\e_k\}$, $\{a_k\}$ and $\{A_k\}$ satisfies \eqref{condition:ekAk}. One may not be able to construct an appropriate initial data by arranging $\{a_1,a_0,a_2,a_0,\cdots\}$ with suitable spatial distance because this may not give a non-BV data if $\sum\limits_{k\geq1}\abs{a_k-a_0}<\f$. To give more oscillation in data and to ensure $TV(u_0)=\f$ we repeat each pair $(a_k,a_0)$ sufficiently many times.  
		For this purpose we define 
		\begin{equation}\label{def:NkJk}
		N_l:=\left[\frac{A_l}{\e_l}\right]\mbox{ for }l\geq1\mbox{ and }J_k:=\sum\limits_{l=1}^{k}N_l\mbox{ for }k\geq1
		\end{equation}
        where $[\cdot]$ is the standard `box function'. Let $J_0=0$. Since $A_l\geq \e_l$ as in \eqref{condition:ekAk} we have $N_l\geq1$. The quantity $N_k$ denotes the number of times the oscillation of strength $\abs{a_k-a_0}$ appears in initial data. Next we consider the following real number sequences $(\si_m)_{m\geq 1}$ and $(\de_m)_{m\geq 1}$ defined as 
		\begin{equation}
		\si_m :=\left\{\begin{array}{lllll}\label{def:sigma}
		a_{0} &\mbox{if}& m=1, &&\\
		a_k &\mbox{if}& m=2i &\mbox{with}& J_{k-1}+1\leq i\leq J_{k},k\geq1,\\
		a_{0} &\mbox{if}& m=2i+1 &\mbox{with}& J_{k-1}+1\leq i\leq J_{k},k\geq1,
        \end{array}\right. 
        \end{equation}
        and 
        \begin{equation}
		\ \de_m :=\left\{\begin{array}{lll}\label{def:delta} 
		T\e_1&\mbox{if}& m=1,2,\\
		T\e_k &\mbox{if}& 2J_{k-1}+3\leq m\leq 2J_{k}+2,\ k\geq1.
		\end{array}\right.
		\end{equation}
		We henceforth write 
		\begin{equation}\label{def:x_n}
		x_m:=\sum\limits_{j=1}^{m}\de_j \mbox{ for } m\geq 1\mbox{ and }x_\f:=\lim\limits_{m\to\f}x_m=\sum\limits_{j=1}^{\f}\de_j.
		\end{equation} 
		Consequently, by using \eqref{condition:ekAk}, $J_k-J_{k-1}=N_k$ and $N_k\e_k\leq2A_k$ we have
		\begin{eqnarray}
		x_\f=\sum\limits_{j=1}^{\f}\de_j&\leq&2T\sum\limits_{k=1}^{\f}(J_{k}-J_{k-1})\e_k+2T\e_1\leq 4T\sum\limits_{k=1}^{\f}A_k+2T\e_1<\f.\label{estimate:x_infty}
		\end{eqnarray}
		We consider an initial data $w_0$ defined as follows  
		\begin{equation}\label{def:initial1}
		w_0(x):=\left\{\begin{array}{rcl}
		a_0&\mbox{ if }&x<x_1,\\
		\si_m&\mbox{ if }&x_{m}<x<x_{m+1},\\
		a_0&\mbox{ if }&x>x_\f .
		\end{array}	\right.
		\end{equation}
		By our assumption $f^{\p\p}>0$ on $(a_0,b_0)$. 
		We observe the following facts: 
		\begin{enumerate}[(a)]
			\item Since $f$ is convex on $(a_0,b_0)$ and $a_k>a_0$, \textit{shock wave} is formed at $(x_{2m+1},0)$ with speed $(f(a_k)-f(a_0))/(a_k-a_0)$ for $J_{k-1}<m\leq J_{k}$ when $k\geq1$.
			
			\item Since $f$ is convex on $(a_0,b_0)$ and $a_k>a_0$ \textit{rarefaction fan} is generated from $(x_{2m},0)$. At small time $t>0$ it covers the interval $\{(x,t);\,x\in(f^{\p}(a_0)t,f^{\p}(a_k)t)\}$ for $J_{k-1}< m\leq J_{k}$ when $k\geq1$.
			
		\end{enumerate}
		We wish to check that no two Riemann problem solutions interact with each other in time interval $[0,T]$.
		\begin{enumerate}[(a)]
			\item\textit{No interaction between shock wave from $(x_{2m+1},0)$ with extreme right characteristic from point $(x_{2m},0)$ for $J_{k-1}<m< J_k$:} Suppose at time $t_1$ shock line, $x=x_{2m+1}+t(f(a_k)-f(a_0))/(a_k-a_0)$ meets with the extreme right characteristic line from $(x_{2m},0)$ i.e. $x=x_{2m}+tf^{\p}(a_k)$. Then we have $x_{2m+1}+t_1(f(a_k)-f(a_0))/(a_k-a_0)=x_{2m}+t_1f^{\p}(a_{k})$. By the Mean Value Theorem we get $t_1=\frac{x_{2m+1}-x_{2m}}{f^{\p}(a_{k})-\frac{f(a_k)-f(a_0)}{a_k-a_0}}=\frac{x_{2m+1}-x_{2m}}{f^{\p}(a_{k})-f^{\p}(c_k)}
			$ for some $c_k\in(a_0,a_k)$. By the choice of $x_{j}$ in \eqref{def:x_n} we have $x_{2m+1}-x_{2m}=\e_k$ for $J_{k-1}<m<J_k$. Since $f^{\p\p}>0$ on $(a_0,b_0)$ and $a_0<c_k<a_k$ we have $f^{\p}(a_0)<f^{\p}(c_k)<f^{\p}(a_k)$. Hence $\abs{f^{\p}(c_k)-f^{\p}(a_k)}=f^{\p}(a_k)-f^{\p}(c_k)<f^{\p}(a_k)-f^{\p}(a_0)$. With the help of \eqref{condition:ekAk} we obtain
			\begin{equation}\label{time_estimate:t1}
			t_1=\frac{x_{2m+1}-x_{2m}}{f^{\p}(a_{k})-f^{\p}(c_k)}\geq\frac{T\e_k}{f^{\p}(a_k)-f^{\p}(a_0)}\geq T.
			\end{equation}
			
			\item\textit{No interaction between shock wave from $(x_{2m+1},0)$ with extreme left characteristic from point $(x_{2m+2},0)$ for $J_{k-1}<m< J_k$:} Suppose at time $t_2$ the shock line $x=x_{2m+1}+t(f(a_k)-f(a_0))/(a_k-a_0)$ meet the extreme left characteristic line from $(x_{2m+2},0)$, i.e. $x=x_{2m+2}+tf^{\p}(a_0)$. Then we have
			$
			x_{2m+1}+t_2\frac{f(a_k)-f(a_0)}{a_k-a_0}=x_{2m+2}+t_2f^{\p}(a_0).
			$ After a rearrangement and by the Mean Value Theorem we have
			$
			t_2=\frac{x_{2m+2}-x_{2m+1}}{\frac{f(a_k)-f(a_0)}{a_k-a_0}-f^{\p}(a_0)}=\frac{x_{2m+2}-x_{2m+1}}{f^{\p}(c_{k})-f^{\p}(a_0)}
			$ for some $c_k\in(a_0,a_k)$. Since $f^{\p\p}>0$ on $(a_0,b_0)$ we have $f^{\p}(a_k)>f^{\p}(c_k)>f^{\p}(a_0)$. Therefore, we obtain $\abs{f^{\p}(c_k)-f^{\p}(a_0)}=f^{\p}(c_k)-f^{\p}(a_0)<f^{\p}(a_k)-f^{\p}(a_0)$. Since $x_{2m+2}-x_{2m+1}=\e_k$ for $J_{k-1}<m<J_k$, by applying \eqref{condition:ekAk}, we have
			\begin{equation}\label{estimate:t2}
			t_2=\frac{x_{2m+2}-x_{2m+1}}{f^{\p}(c_{k})-f^{\p}(a_0)}\geq\frac{T\e_k}{f^{\p}(a_k)-f^{\p}(a_0)}\geq T.
			\end{equation}
			
			\item\textit{No interaction between shock wave from $(x_{2m+1},0)$ with extreme right characteristic from point $(x_{2m},0)$ for $m= J_k$:} Suppose $x=x_{2m+1}+t(f(a_k)-f(a_0))/(a_k-a_0)$ meets with the line $x=x_{2m}+tf^{\p}(a_k)$ at time $t_3$. By similar argument as in \eqref{time_estimate:t1}, we have $t_3\geq T$.
			
			\item\textit{No interaction between shock wave from $(x_{2m+1},0)$ with extreme left characteristic from point $(x_{2m+2},0)$ for $m= J_k$:} Suppose $x=x_{2m+1}+t(f(a_k)-f(a_0))/(a_k-a_0)$ meets with the line $x=x_{2m+2}+tf^{\p}(a_0)$ at time $t_4$.  
			By our choice of $x_{2m+2}$ for $m=J_k$ we have $x_{2m+2}-x_{2m+1}=\e_k$. Now by using \eqref{condition:ekAk} with a similar argument as in \eqref{estimate:t2} we infer that $t_4\geq T$.
%
%
		\end{enumerate}
		
				\begin{figure}
			\centering
			\begin{tikzpicture}[scale=.7]

			\draw[thick,->] (-8 ,0) -- (7,0) node[anchor=north west] {$x$};
			\draw[thick,->] (-8 ,0) -- (-8,6.5) node[anchor=north east] {$t$};
			
			\draw[color=black] (-7.2 ,0) -- (-6.2,5) ;
			\draw[color=black] (-7.2 ,0) -- (-6.1,5) ;
			\draw[color=black] (-7.2 ,0) -- (-6.0,5) ;
			\draw[color=black] (-7.2 ,0) -- (-5.9,5) ;
			\draw[color=black] (-7.2 ,0) -- (-5.8,5) ;
			\draw[color=black] (-7.2 ,0) -- (-5.7,5) ;
			\draw[color=black] (-7.2 ,0) -- (-5.6,5) ;
			\draw[color=black] (-7.2 ,0) -- (-5.5,5) ;
			\draw[color=black] (-7.2 ,0) -- (-5.4,5) ;
			\draw[color=black] (-7.2 ,0) -- (-5.3,5) ;
			\draw[color=black] (-7.2 ,0) -- (-5.2,5) ;
			\draw[color=black] (-7.2 ,0) -- (-5.1,5) ;
			\draw[color=black] (-7.2 ,0) -- (-5.0,5) ;
			\draw[color=black] (-7.2 ,0) -- (-4.9,5) ;
			\draw[color=black] (-7.2 ,0) -- (-4.8,5) ;
			\draw[color=blue] (-6.5 ,0) -- (-4.8,5) ;
			\draw[color=black] (-5.8 ,0) -- (-4.8,5) ;
			\draw[color=black] (-5.8 ,0) -- (-4.7,5) ;
			\draw[color=black] (-5.8 ,0) -- (-4.6,5) ;
			\draw[color=black] (-5.8 ,0) -- (-4.5,5) ;
			\draw[color=black] (-5.8 ,0) -- (-4.4,5) ;
			\draw[color=black] (-5.8 ,0) -- (-4.3,5) ;
			\draw[color=black] (-5.8 ,0) -- (-4.2,5) ;
			\draw[color=black] (-5.8 ,0) -- (-4.1,5) ;
			\draw[color=black] (-5.8 ,0) -- (-4.0,5) ;
			\draw[color=black] (-5.8 ,0) -- (-3.9,5) ;
			\draw[color=black] (-5.8 ,0) -- (-3.8,5) ;
			\draw[color=black] (-5.8 ,0) -- (-3.7,5) ;
			\draw[color=black] (-5.8 ,0) -- (-3.6,5) ;
			\draw[color=black] (-5.8 ,0) -- (-3.5,5) ;
			\draw[color=black] (-5.8 ,0) -- (-3.4,5) ;
			\draw[color=blue] (-5.1 ,0) -- (-3.4,5) ;
			\draw[color=black] (-4.4 ,0) -- (-3.4,5) ;
			\draw[color=black] (-3.2 ,0) -- (-2.2,5) ;
			\draw[color=black] (-3.2 ,0) -- (-2.1,5) ;
			\draw[color=black] (-3.2 ,0) -- (-1.976,5.1) ;
			\draw[color=black] (-3.2 ,0) -- (-1.848,5.2) ;
			\draw[color=black] (-3.2 ,0) -- (-1.716,5.3) ;
			\draw[color=black] (-3.2 ,0) -- (-1.58,5.4) ;
			\draw[color=black] (-3.2 ,0) -- (-1.44,5.5) ;
			\draw[color=black] (-3.2 ,0) -- (-1.296,5.6) ;
			\draw[color=black] (-3.2 ,0) -- (-1.148,5.7) ;
			\draw[color=black] (-3.2 ,0) -- (-0.996,5.8) ;
			\draw[color=black] (-3.2 ,0) -- (-0.84,5.9) ;
			\draw[color=black] (-3.2 ,0) -- (-0.68,6) ;
			\draw[color=black] (-3.2 ,0) -- (-0.56,6) ;
			\draw[color=black] (-3.2 ,0) -- (-0.44,6) ;
			\draw[color=black] (-3.2 ,0) -- (-0.8,5) ;
			\draw[color=blue] (-2.5 ,0) -- (-0.8,5) ;
			\draw[color=black] (-1.8 ,0) -- (-0.8,5) ;

			\draw[color=black] (-0.8 ,0) -- (0,6) ;
\draw[color=black] (-0.8 ,0) -- (0.1,6) ;
\draw[color=black] (-0.8 ,0) -- (0.2,6) ;
\draw[color=black] (-0.8 ,0) -- (0.3,6) ;
\draw[color=black] (-0.8 ,0) -- (0.4,6) ;
\draw[color=black] (-0.8 ,0) -- (0.5,6) ;
\draw[color=black] (-0.8 ,0) -- (0.6,6) ;
\draw[color=black] (-0.8 ,0) -- (0.7,6) ;
\draw[color=black] (-0.8 ,0) -- (0.8,6) ;
\draw[color=blue] (-0.4 ,0) -- (0.8,6) ;
\draw[color=black] (0 ,0) -- (0.8,6) ;
\draw[color=black] (0 ,0) -- (0.9,6) ;
\draw[color=black] (0 ,0) -- (1.0,6) ;
\draw[color=black] (0 ,0) -- (1.1,6) ;
\draw[color=black] (0 ,0) -- (1.2,6) ;
\draw[color=black] (0 ,0) -- (1.3,6) ;
\draw[color=black] (0 ,0) -- (1.4,6) ;
\draw[color=black] (0 ,0) -- (1.5,6) ;
\draw[color=black] (0 ,0) -- (1.6,6) ;
\draw[color=blue] (0.4 ,0) -- (1.6,6) ;
\draw[color=black] (0.8 ,0) -- (1.6,6) ;
\draw[color=black] (2.2 ,0) -- (3,6) ;
\draw[color=black] (2.2 ,0) -- (3.1,6) ;
\draw[color=black] (2.2 ,0) -- (3.2,6) ;
\draw[color=black] (2.2 ,0) -- (3.3,6) ;
\draw[color=black] (2.2 ,0) -- (3.4,6) ;

\draw[color=black] (2.2 ,0) -- (3.5,6) ;
\draw[color=black] (2.2 ,0) -- (3.6,6) ;
\draw[color=black] (2.2 ,0) -- (3.7,6) ;
\draw[color=black] (2.2 ,0) -- (3.8,6) ;
\draw[color=blue] (2.6 ,0) -- (3.8,6) ;
\draw[color=black] (3 ,0) -- (3.8,6) ;
			
			\draw[color=blue] plot [smooth] coordinates { (-0.8,5) (-0.6,5.2) (-0.4,5.7) (-0.2,6)};

			\draw[dashed] (6,4.5) -- (-8,4.5) node[anchor= east] {$T$};
			
			\draw[dashed,<->] (-6.5,-0.6) -- (-2.5,-0.6) ;
			\draw[thick][] (-4.5,-.6) node[anchor=north] {$N_1$ times};
			\draw[dashed,<->] (-0.4,-0.6) -- (2.6,-0.6) ;
			\draw[thick][] (1.1,-0.6) node[anchor=north] {$N_2$ times};
			\draw[thick][] (-7.5,0.52) node[anchor=north] {$a_0$};
			
			\draw[thick][] (-6.7,0.52) node[anchor=north] {$a_1$};
			\draw[thick][] (-6.05,0.52) node[anchor=north] {$a_0$};
			
			\draw[thick][] (-5.4,0.52) node[anchor=north] {$a_1$};
			\draw[thick][] (-4.75,0.52) node[anchor=north] {$a_0$};
			
			\draw[thick][] (-2.8,0.52) node[anchor=north] {$a_1$};
			\draw[thick][] (-2.05,0.52) node[anchor=north] {$a_0$};
			
			\draw[thick][] (-1.3,0.52) node[anchor=north] {$a_0$};
			
			\draw[thick][] (-.9,0.52) node[anchor=north,  transform canvas={scale=.65}] {$a_2$};
			\draw[thick][] (-0.3,0.52) node[anchor=north,  transform canvas={scale=.65}] {$a_0$};
			
			\draw[thick][] (0.36,0.52) node[anchor=north,  transform canvas={scale=.65}] {$a_2$};
			\draw[thick][] (0.98,0.52) node[anchor=north,  transform canvas={scale=.65}] {$a_0$};
			
			\draw[thick][] (3.72,0.52) node[anchor=north,  transform canvas={scale=.65}] {$a_2$};
			\draw[thick][] (4.35,0.52) node[anchor=north,  transform canvas={scale=.65}] {$a_0$};
			
			\draw[thick][] (-2.6,1.4) node[anchor=south west, transform canvas={scale=1.5}]{$\textbf{...} $};
			\draw[thick][] (0.84,1.4) node[anchor=south west, transform canvas={scale=1.5}]{$\textbf{...} $};
			\draw[thick][] (2.7,1.4) node[anchor=south west, transform canvas={scale=1.5}]{$\textbf{...} $};
			
			\draw[thick][] (3.7,-0.23) node[anchor=south west, transform canvas={scale=1.5}]{$\textbf{.} $};
			
			\draw[thick][] (1.48,-0.23) node[anchor=south west, transform canvas={scale=1.5}]{$\textbf{.} $};
			\draw[thick][] (-1.93,-0.23) node[anchor=south west, transform canvas={scale=1.5}]{$\textbf{.} $};
			\draw[thick][] (-0.79,-0.23) node[anchor=south west, transform canvas={scale=1.5}]{$\textbf{.} $};
			\draw[thick][] (-5.06,-0.23) node[anchor=south west, transform canvas={scale=1.5}]{$\textbf{.} $};
			
			\draw[thick][] (-2.5,0) node[anchor=north] {$x_{2J_1+1}$};
			\draw[thick][] (-0.8,0) node[anchor=north] {$x_{2J_1+2}$};
			\draw[thick][] (2.6,0) node[anchor=north] {$x_{2J_2+1}$};
			\draw[thick][] (-7.2,0) node[anchor=north] {$x_{2}$};
			\draw[thick][] (6.1,-.1) node[anchor=north] {$x_\f$};
			\end{tikzpicture}
			\caption{The initial data is constructed by repeating the pair $(a_0,a_k)$ $N_k$ times and two consecutive pair are connected by $a_0$. At the $(a_0,a_k)$ jump, rarefaction arises in structure of the associated entropy solution. Shocks are formed at each jump point $(a_k,a_0)$. Waves do not interact till time $T$.}\label{figure:nonlinear}
		\end{figure}
		Above observation confirms that no two Riemann problem solutions interact with each other before time $T$ (see Figure \ref{figure:nonlinear} for clear illustration). Let $w$ be the entropy solution to \eqref{eqn:conlaw} with initial data $w_0$. Then we have
		\begin{equation}\label{structure:1Dsol}
		w(x,t)=\left\{\begin{array}{lll}
		a_0&\mbox{if}&\ x<x_2+f^{\p}(a_0)t,\\
		a_{k} &\mbox{if}& \begin{array}{lll} x_{2i}+f^{\p}(a_{k})t<x<x_{2i+1}+P_kt,\\
		\mbox{for }  J_{k-1}+1\leq i\leq J_{k},\ k\geq 1,
		\end{array}\\
		a_0 &\mbox{if}&\begin{array}{lll} x_{2i+1}+P_kt<x<x_{2i+2}+f^{\p}(a_0)t,\\
		\mbox{for }  J_{k-1}+1\leq i\leq J_{k},\ k\geq 1, 
		\end{array}\\
		(f^{\p})^{-1}\left(\frac{x-x_{2i}}{t}\right) &\mbox{if}& \begin{array}{lll} x_{2i}+f^{\p}(a_0)t<x\leq x_{2i}+f^{\p}(a_{k})t,\\
		\mbox{for }  J_{k-1}+1\leq i\leq J_{k},\ k\geq 1,
		\end{array}\\
		a_0&\mbox{if}&\ x>x_\f+f^{\p}(a_0)t,
		\end{array}\right.
		\end{equation}
		for $t\in[0,T]$, where $P_k=(f(a_k)-f(a_0))/(a_k-a_0)$.	Now we wish to show that $TV(w(\cdot,t),[(x_1+f^{\p}(a_0)t,x_\f+f^{\p}(a_0)t)])=\f$ for all $t\in[0,T]$. To do this, fix a time $t\in[0,T]$. Consider a sequence $\{z_n\}_{n\geq1}$ such that
			\begin{equation*}
		z_n\in\left\{\begin{array}{lllll}
		(x_{2i}+f^{\p}(a_{k})t,x_{2i+1}+P_kt)&\mbox{if}& n=2i+1 &\mbox{with}& J_{k-1}<i\leq J_{k},\\
		&&&&\mbox{ for }k\geq1,\\
		(x_{2i+1}+P_kt,x_{2i+2}+f^{\p}(a_0)t) &\mbox{if}&n=2i+2 &\mbox{with}& J_{k-1}< i\leq J_{k},\\
		&&&&\mbox{ for }k\geq1.\\
		\end{array}\right.
		\end{equation*}
        Note that $\{z_n\}\subset[x_1+f^{\p}(a_0)t,x_\f+f^{\p}(a_0)t]$. From \eqref{structure:1Dsol} we have
		\begin{equation}\label{nonlinear:w(z_n,t)}
		w(z_n,t)\in\left\{\begin{array}{lllll}
		a_k&\mbox{if}& n=2i+1 &\mbox{with}& J_{k-1}<i\leq J_{k},\\
		a_0&\mbox{if}&n=2i+2 &\mbox{with}& J_{k-1}< i\leq J_{k},
		\end{array}\right.\mbox{ for }k\geq1,
		\end{equation}
	
		therefore, from \eqref{nonlinear:w(z_n,t)} we infer
		\begin{align}
		\sum\limits_{j=3}^{\f}\abs{w(z_j,t)-w(z_{j+1},t)}
		&=\sum\limits_{k=1}^{\f}\sum\limits_{i=J_{k-1}+1}^{J_k}\abs{w(z_{2i+2},t)-w(z_{2i+1},t)}\nonumber\\
		&=\sum\limits_{k=1}^{\f}\sum\limits_{i=J_{k-1}+1}^{J_k}\abs{a_k-a_0}=\sum\limits_{k=1}^{\f}(J_k-J_{k-1})\abs{a_k-a_0}.\label{nonlinear:cal1}
		\end{align}
		Recall from \eqref{def:NkJk} that $J_k-J_{k-1}=N_k$ and $N_k$ is the largest integer smaller than or equal to $\e_k^{-1}\abs{f^{\p}(a_k)-f^{\p}(a_0)}/\abs{a_k-a_0}$. Therefore, $N_k\geq \frac{\abs{f^{\p}(a_k)-f^{\p}(a_0)}}{2\e_k\abs{a_k-a_0}}$. By using the above observation we can estimate \eqref{nonlinear:cal1} as follows
		\begin{align}
		\sum\limits_{j=1}^{\f}\abs{w(z_j,t)-w(z_{j+1},t)}
		\geq\sum\limits_{k=1}^{\f}\frac{\abs{f^{\p}(a_k)-f^{\p}(a_0)}}{2\e_k\abs{a_k-a_0}}\abs{a_k-a_0}
		&\geq \sum\limits_{k=1}^{\f}\frac{\abs{f^{\p}(a_k)-f^{\p}(a_0)}}{2\e_k}.\label{nonlinear:cal2}
		\end{align}
		Since $\e_k=\abs{a_k-a_0}^{p_k}$ we get 
		\begin{equation*}
		(2\e_k)^{-1}\abs{f^{\p}(a_k)-f^{\p}(a_0)}=(2\abs{a_k-a_0}^{p_k})^{-1}\abs{f^{\p}(a_k)-f^{\p}(a_0)}.
		\end{equation*}From \eqref{eq:sequence0} we obtain $(\abs{a_k-a_0}^{p_k})^{-1}\abs{f^{\p}(a_k)-f^{\p}(a_0)}=k^{-1}$. Therefore, from \eqref{nonlinear:cal2} we have
		\begin{equation*}
		\sum\limits_{j=1}^{\f}\abs{w(z_j,t)-w(z_{j+1},t)}
		\geq\sum\limits_{k=1}^{\f}\frac{\abs{f^{\p}(a_k)-f^{\p}(a_0)}}{2\abs{a_k-a_0}^{p_k}}=\sum\limits_{k=1}^{\f}\frac{1}{2k}=\f.
		\end{equation*}
		This shows that \begin{equation}\label{TVblowup:w1}
		TV(w(\cdot,t),[(x_1+f^{\p}(a_0)t,x_\f+f^{\p}(a_0)t)])=\f\mbox{ for }t\in[0,T]. 
		\end{equation}
		To conclude Proposition \ref{Proposition:1d_uptoT} we need a compactly supported initial data. Note that $w_0$ as in \eqref{def:initial1} could be such candidate if $a_0=0$. In general, $a_0$ may not be $0$. To obtain a compactly supported data, we apply similar technique as we have done in \descref{linear}{Case(1)}. Let $M_2$ be defined as follows
		\begin{equation}\label{def:M2}
		M_2:=\max\{\abs{f^{\p}(u)};\,\abs{u}\leq r_2\}\mbox{ where }r_2:=\|w_0\|_{L^{\f}(\R)}.
		\end{equation}
		From \eqref{def:initial1} observe that $\|w_0\|_{L^{\f}(\R)}=\max\{\abs{a_0},\abs{a_k};\mbox{ for }k\geq1\}$. Therefore $r_2$ does not depend on $T$. Hence $M_2$ is independent of $T$. Let $R_2=3M_2T$ and $I_1=[x_1-R_2,x_\f+R_2]$. We define 
		\begin{equation}\label{def:v_0_nonlin}
		v_0^T:=w_0\chi_{I_1}
		\end{equation}
		where $w_0$ is defined as in \eqref{def:initial1}. Note that by choice of $x_k$ in \eqref{def:x_n} we have $x_1\geq0$. Recall the constant $C$ defined as in \eqref{condition:ekAk}. From \eqref{estimate:x_infty} we have $x_\f\leq 4TC$. Hence we have 
		\begin{equation}\label{eq:support_nonlin}
		supp(v_0^T)\subset [-3M_2T,(C_1+3M_2)T]\mbox{ where }C_1:=4C.
		\end{equation}Let $v^T$ be the entropy solution to \eqref{eqn:conlaw} with initial data $v_0^T$. Applying Proposition \ref{Prop:Kruzkov} we have for $t\in[0,T]$ and $z_1<z_2$
		\begin{equation}\label{nonlinear_case1}
		\int\limits_{[z_1,z_2]}\abs{w(x,t)-v^T(x,t)}\,dx\leq \int\limits_{[z_1-M_2T,z_2+M_2T]}\abs{w_0(x)-v^T_0(x)}\,dx.
		\end{equation}
		Set $z_1=x_1-M_2T$ and $z_2=x_\f+M_2T$. With this choice of $z_1$ and $z_2$ we observe that $[z_1-M_2T,z_2+M_2T]\subset I_1$. Since $v_0^T=w_0$ in $I_1$. From \eqref{nonlinear_case1} we have 
		\begin{equation}\label{nonlinear:w=v}
		w(x,t)=v^T(x,t) \mbox{ for all }(x,t)\in [x_1-M_2T,x_\f+M_2T]\times[0,T]. 
		\end{equation}
		From \eqref{def:M2} we have $M_2\geq \abs{f^{\p}(a_0)}$. Therefore, we have
		\begin{equation}\label{nonlinear:cal3}
		[x_1+f^{\p}(a_0)t,x_\f+f^{\p}(a_0)t]\subset[x_1-M_2T,x_\f+M_2T]\mbox{ for all }t\in[0,T].
		\end{equation}
		By \eqref{TVblowup:w1}, \eqref{nonlinear:w=v} and \eqref{nonlinear:cal3} we have 
		$TV(v^T(\cdot,t),[x_1-M_2T,x_\f+M_2T])=\f$ for all $ t\in[0,T]$. Recall the definition of $C_1$ in \eqref{eq:support_nonlin} and $x_\f\leq 2CT$. Hence $x_\f+M_2T\leq (C_1+M_2)T$. Since $x_1>0$ we have $[x_1-M_2T,x_\f+M_2T]\subset[-3M_2T,(C_1+3M_2)T]$. Therefore, we obtain
		\begin{equation}\label{TVblowup_nonlin}
		TV(v^T(\cdot,t),[-3M_2T,(C_1+3M_2)T])=\f\mbox{ for all } t\in[0,T].
		\end{equation}
	\end{description}
	Set $r=\max\{r_1,r_2\}$ where $r_1,r_2$ are defined as in \eqref{def:M1} and \eqref{def:M2} respectively. Let $v_{0,T}^{lin},v_0^T$ as in \eqref{def:v_0lin} and \eqref{def:v_0_nonlin} respectively. Subsequently, $\|v_{0,T}^{lin}\|_{L^{\f}(\R)},\|v^T_0\|_{L^{\f}(\R)}\leq r$. Recall that $v_{lin}^T$ is the entropy solution to \eqref{eqn:conlaw} with initial data $v_{0,T}^{lin}$ as in \descref{linear}{Case(1)} and $v^T$ is the entropy solution to \eqref{eqn:conlaw} with initial data $v_{0}^{T}$ as in \descref{non-linear}{Case(2)}. By maximum principle $\|v_{lin}^T\|_{L^{\f}(\R)},\|v^T\|_{L^{\f}(\R)}\leq r$. As in \descref{linear}{Case(1)} and \descref{non-linear}{Case(2)} note that $r_1$ and $r_2$ are independent of $T$ and so is $r$. Let us define $M:=\max\{\abs{f^{\p}(u)};\,\abs{u}\leq r\}$. Since $r\geq r_1,r_2$ we have $M\geq M_1, M_2$. Observe that $M$ is independent of $T$ because $r$ does not depend on $T$. Note that $y_\f$ in \descref{linear}{Case(1)} and $C_1$ in \eqref{eq:support_nonlin} are independent of $T$. Consider $\la:=y_\f+3M+C_1$ which is independent of $T$. Note that $[-3M_1T,y_\f+3M_1T]\subset[-\la T,\la(T+1)]$ and $[-3M_2T,(C_1+3M_2)T]\subset[-\la T,\la(T+1)]$. Combining \eqref{eq:support_lin} and \eqref{eq:support_nonlin} we can write $supp(v_{0,T}^{lin}),\,supp(v^T_0)\subset [-\la T,\la(T+1)]$.
	From \eqref{TVblowup_lin} and \eqref{TVblowup_nonlin} we infer
	$ 
	TV(v_{lin}^T(\cdot,t),[-\la T,\la(T+1)]),TV(v^T(\cdot,t),[-\la T,\la(T+1)])=\f\mbox{ for all }t\in[0,T].
	$ 
	This concludes Proposition \ref{Proposition:1d_uptoT}.
\end{proof}

\subsubsection{Proof of Theorem \ref{theorem:1d}}
Now we are ready to prove Theorem \ref{theorem:1d}. From Proposition \ref{Proposition:1d_uptoT} there exists an entropy solution such that $TV(v^T(\cdot,t))=\f$ for $0\leq t<T$. By using finite speed of propagation we construct a solution $v$ which does not belong to $BV_{loc}(\R)$ for all time $t>0$.
\begin{proof}[Proof of Theorem \ref{theorem:1d}]
	Since there exists an $a_0\in\R$ such that $f^{\p\p}(a_0)=0$ we can invoke Proposition \ref{Proposition:1d_uptoT}. Consider a positive increasing sequence $\{T_n\}$ such that $T_n\rr\f$ as $n\rr\f$. By using Proposition \ref{Proposition:1d_uptoT} there exists $\la>0$ and $u^*$ such that for each $T_n$ we have an initial data $v_0^{T_n}$ such that it satisfies the following:
	\begin{enumerate}
		\item $\|v^{T_n}_0\|_{L^{\f}(\R)}\leq u^*$ and $supp(v^{T_n}_0)\subset[-\la T_n,\la(T_n+1)]$.
		\item if $v^{T_n}$ is the entropy solution to \eqref{eqn:conlaw} with initial data $v^{T_n}_0$, then 
		\begin{equation}\label{TVblow_up1}
		TV(v^{T_n}(\cdot,t),[-\la T_n,\la(T_n+1)])=\f\mbox{ for all }t\in[0,T_n].
		\end{equation}
	\end{enumerate}
We want to build an entropy solution so that it remains non-$BV_{loc}(\R)$ for all time. To do this, we use the sequence of data $v^{T_n}_0$ and with the help of finite speed of propagation, we arrange them such a way that up to time $T_n$ the entropy solution $v^{T_n}$ does not interact with $v^{T_{m}}$ if $m\neq n$. To this end, let us define $M$ as follows $
	M:=\max\{\abs{f^{\p}(w)};\,\abs{w}\leq u^*\}$. Now define $\nu_n$ as follows 
	\begin{equation}\label{def:nu_n}
	\nu_n=\sum\limits_{k=1}^{n}4(T_k+1)(\la+M)+2(T_{n+1}+1)(\la+M).
	\end{equation}
	Consider the following initial data $v_0$ and $V_0^m$
	\begin{equation}\label{def:initial1D_all_t}
	v_0(x)=\sum\limits_{n=1}^{\f}v_0^{T_n}(x-\nu_n)\mbox{ and }V_0^m(x)=\sum\limits_{n=1}^{m}v_0^{T_n}(x-\nu_n).
	\end{equation}
	Let $v$ and $V^m$ be the entropy solutions to \eqref{eqn:conlaw} corresponding to $v_0$ and $V_0^m$ respectively. To prove $v(\cdot,t)\notin BV_{loc}(\R)$ for all $t>0$, we first show that for $t\in[0,T_{m+1}]$ $supp(v^{T_{m+1}}(\cdot-\nu_{m+1},t))$ is disjoint from $supp(V^m(\cdot,t))$ and $supp(v^{T_{m+k}}(\cdot-\nu_{m+k},t))$ with $k\geq2$. Then we prove that $v(\cdot,t)$ matches with $v^{T_{m+1}}(\cdot-\nu_{m+1},t)$ on some interval for $t\in[0,T_{m+1}]$ and finally by \eqref{TVblow_up1} we conclude the result. Our next aim is to show that $v^{T_{m+1}}(\cdot-\nu_{m+1},t)$ and $V^m(\cdot,t)$ have disjoint support up to time $T_{m+1}$. From \eqref{def:nu_n},
	\begin{equation}\label{difference:nu}
	\nu_{n}-\nu_{n-1}=2(T_n+1)(\la+M)+2(T_{n+1}+1)(\la+M).
	\end{equation}
	Hence
	\begin{align}
	\nu_{n-1}+\la(T_{n-1}+1)
	&\leq \nu_{n-1}+\la(T_n+1)\nonumber\\
	&=-(\nu_{n}-\nu_{n-1})+\la(T_n+1)+\nu_n\nonumber\\
	&=-2(T_n+1)(\la+M)-2(T_{n+1}+1)(\la+M)+\la(T_n+1)+\nu_n\nonumber\\
	&=-\la T_n+\nu_n-2(\la+M(T_n+1)+2(T_{n+1}+1)(\la+M))\nonumber\\
	&\leq-\la T_n+\nu_n.\label{1D:cal1}
	\end{align}Since $supp(v_0^{T_n})$ lies in $[-\la T_n,\la(T_n+1)]$ we have 
	\begin{equation}\label{supp:transalation}
	supp(v_0^{T_n}(\cdot-\nu_n))\subset[-\la T_n+\nu_n,\la(T_n+1)+\nu_n].
	\end{equation}
	From \eqref{def:initial1D_all_t}, \eqref{1D:cal1} and \eqref{supp:transalation} we infer
	\begin{equation*}
	supp(V^m_0)\subset \bigcup\limits_{n=1}^{m}[-\la T_n+\nu_n,\la(T_n+1)+\nu_n]\subset [-\la_1 T_1+\nu_1,\la(T_{m}+1)+\nu_m].
	\end{equation*}
	By \eqref{1D:cal1} and \eqref{supp:transalation}, we get $supp(v_0^{T_n})$ and $supp(v_0^{T_m})$ are disjoint if $n\neq m$. By using the previous observation, from \eqref{def:initial1D_all_t} we observe that $\|v_0\|_{L^{\f}(\R)},\|V^m_0\|_{L^{\f}(\R)}\leq u^*$ for all $m\geq1$ since $\|v^{T_n}_0\|_{L^{\f}(\R)}\leq u^*$  for all $n\geq1$. Therefore, by Proposition \ref{Prop:Kruzkov} the entropy solution $V^m$ to \eqref{eqn:conlaw} with data $V_0^m$ satisfies
	\begin{equation}\label{support:soln_V}
	supp(V^m(\cdot,t))\subset [-\la T_1+\nu_1-MT_{m+1},\la(T_m+1)+\nu_m+MT_{m+1}],
	\end{equation}
	$\mbox{ for all }t\in[0,T_{m+1}]$. Recall that $v^{T_{m+1}}$ is the entropy solution to \eqref{eqn:conlaw} with initial data $v_0^{T_{m+1}}$. By Lemma \ref{lemma:linear_transform}, we infer that $v^{T_{m+1}}(\cdot-\nu_{m+1},t)$ is the entropy solution to \eqref{eqn:conlaw} associated with initial data $v_0^{T_{m+1}}(\cdot-\nu_{m+1})$. As we observed in \eqref{supp:transalation} $v_0^{T_{m+1}}(\cdot-\nu_{m+1})$ has support in $[-\la T_{m+1}+\nu_{m+1},\la(T_{m+1}+1)+\nu_{m+1}]$. Again by Lemma \ref{lemma:linear_transform} we have
	\begin{equation}\label{support:vT}
	supp(v^{T_{m+1}}(\cdot-\nu_{m+1},t))\subset [A_{m+1},B_{m+1}]\mbox{ for all }t\in[0,T_{m+1}]
	\end{equation}
	where $A_{m+1},B_{m+1}$ are defined as follows
	\begin{equation}\label{def:AmBm}
	A_{m+1}=-\la T_{m+1}+\nu_{m+1}-MT_{m+1},\,B_{m+1}=\la(T_{m+1}+1)+\nu_{m+1}+MT_{m+1}.
	\end{equation} 
	From definition \eqref{def:nu_n} of $\nu_n$ and \eqref{difference:nu}, 
	\begin{align*}
	-\la T_{m+1}+\nu_{m+1}-MT_{m+1}
	&=-(\la+M)T_{m+1}+\nu_{m+1}-\nu_m+\nu_m\nonumber\\
	&=-(\la+M)T_{m+1}+2(T_{m+1}+1)(\la+M)\nonumber\\
	&+(T_{m+2}+1)(\la+M)+\nu_m\nonumber\\
	&=(T_{m+1}+2)(\la+M)+(T_{m+2}+1)(\la+M)+\nu_m\nonumber\\
	&=\nu_m+MT_{m+1}+\la(T_{m+2}+1)+\la(T_{m+1}+2)\nonumber\\
	&+M(T_{m+2}+3).
	\end{align*}
	Since $\{T_m\}$ is increasing and $\la,M>0$, 
	\begin{equation}\label{inequality:cal1}
	-\la T_{m+1}+\nu_{m+1}-MT_{m+1}>\la(T_m+1)+\nu_m+MT_{m+1}.
	\end{equation}
	Combining \eqref{support:soln_V}, \eqref{support:vT} and \eqref{inequality:cal1} we obtain
	\begin{equation}\label{disjoint_support_pre}
	supp(V^m(\cdot,t))\cap supp(v^{T_{m+1}}(\cdot-\nu_{m+1},t))=\emptyset\mbox{ for all }t\in[0,T_{m+1}].
	\end{equation}
	Next we show that $v^{T_{m+1}}(\cdot-\nu_{m+1},t)$ and $v^{T_{m+k}}(\cdot-\nu_{m+k},t)$ have disjoint support for $t\in[0,T_{m+1}]$ and $k\geq2$. To this end we consider $m+k$ instead of $m$ in \eqref{disjoint_support_pre} to conclude that $supp(V^{m+k}(\cdot,t))$ and $supp(v^{T_{m+k+1}}(\cdot-\nu_{m+k+1},t))$ are disjoint for $t\in[0,T_{m+k+1}]$. Since $\{T_n\}$ is increasing sequence, we have $T_{m+k+1}\geq T_{m+1}$ for all $k\geq1$. By previous observation we have $supp(v^{T_{m+k+1}}(\cdot-\nu_{m+k},t))$ and $supp(V^{m+k}(\cdot,t))$ are disjoint for $t\in[0,T_{m+1}]$ and $k\geq1$. We also know that support of $v^{T_{m+1}}(\cdot-\nu_{m+1},t)$ is contained in support of $V^{m+k}(\cdot,t)$ for $k\geq1$. Hence 
	\begin{equation}\label{disjoint_support_post}
	supp(v^{T_{m+k}}(\cdot-\nu_{m+k},t))\cap supp(v^{T_{m+1}}(\cdot-\nu_{m+1},t))=\emptyset
	\end{equation}
	$\mbox{ for all }t\in[0,T_{m+1}]\mbox{ and }k\geq2$. From \eqref{disjoint_support_pre} and \eqref{disjoint_support_post} we conclude that $v^{T_{m+1}}(\cdot-\nu_{m+1},t)$ does not interact with $V^{m}$ and $v^{T_{m+k}}(\cdot-\nu_{m+k},t)$ for $t\in[0,T_{m+1}]$ and $k\geq2$. Recall that $v$ is the entropy solution to \eqref{eqn:conlaw} with initial data $v_0$ where $v_0$ is defined as in \eqref{def:initial1D_all_t}. Therefore,
	$
	v(x,t)=v^{T_{m+1}}(x-\nu_{m+1},t)$ for $(x,t)\in[A_{m+1},B_{m+1}]\times[0,T_{m+1}]$. 
	 By using \eqref{TVblow_up1} and \eqref{def:AmBm} we have for $t\in[0,T_{m+1}]$, 
	$$
	TV(v(\cdot,t),[A_{m+1},B_{m+1}])\geq TV(v^{T_{m+1}}(\cdot-\nu_{m+1},t),[-\la T_{m+1},\la(T_{m+1}+1)])=\f.
	$$
	This concludes Theorem \ref{theorem:1d}.
\end{proof}	

\subsection{Construction of multi-D non-BV solution}
In this section, we focus on proof of Theorem \ref{theorem1}. It follows in the same way as proof of Theorem \ref{theorem:1d}. In the proof of Theorem \ref{theorem:1d}, we have used a compact support data $v_0^T$ as a building block which is guaranteed by Proposition \ref{Proposition:1d_uptoT}. Similarly, in multi dimension we first construct a compact support data $U_{\xi,T}^0$ such that the associated entropy solution $U_{\xi,T}$ does not belong to $BV_{loc}(\R^d)$ for $t\in[0,T]$. By using $U_{\xi,T}^0$ as a building block, we prove Theorem \ref{theorem1} in the same manner as Theorem \ref{theorem:1d}.
\begin{proof}[Proof of Theorem \ref{theorem1}]
	Fix $\hat{b}\in\re$. Now we wish to find a direction $\xi\in\mathcal{S}^{d-1}$ such that $F^{\p\p}(\hat{b})\cdot\xi=0$. If $F^{\p\p}(\hat{b})=\textbf{0}\in\re^d$ then we can take $\xi$ to be any vector in $\mathcal{S}^{d-1}$. If $F^{\p\p}(\hat{b})\neq\textbf{0}$ then there exists $i\in\{1,\cdots,d\}$ such that $F^{\p\p}_i(a_0)\neq0$. Therefore, in this case we choose $\xi=(\xi_{1},\cdots,\xi_{d})\in\mathcal{S}^{d-1}$ as follows:
	\begin{equation*}
	\xi_{j}:=\left\{\begin{array}{lll}
	-\vartheta^{-1/2} F^{\p\p}_i(\hat{b})&\mbox{for}&j\neq i,\\
	\vartheta^{-1/2}\sum\limits_{j\neq i}F^{\p\p}_j(\hat{b})&\mbox{for}&j=i,
	\end{array}\right.
	\mbox{ with }
	\vartheta=(d-1)\abs{F^{\p\p}_i(\hat{b})}^2+\abs{\sum\limits_{j\neq i}F^{\p\p}_j(\hat{b})}^2.
	\end{equation*}
	Consider $f=F\cdot\xi$. Then $f^{\p\p}(\hat{b})=0$. Since $F\in C^2(\R,\R^d)$ we have $f\in C^2(\R,\R)$. Similar to Theorem \ref{theorem:1d} consider a positive increasing real number sequence $\{T_n\}$ such that $T_n\rr\f$ as $n\rr\f$. Now invoke Proposition \ref{Proposition:1d_uptoT} to get generic constant $\la,u^*>0$ such that for each $T_n$ we get an initial data $v^{T_n}_0$ satisfying the following properties:
	\begin{enumerate}
		\item $\|v^{T_n}_0\|_{L^{\f}(\R)}\leq u^*$ and $
		supp(v^{T_n}_0)\subset[-\la T_n,\la(T_n+1)] $.
		\item if $v^{T_n}$ is the entropy solution to \eqref{eqn:conlaw} with initial data $v^{T_n}_0$, then 
		\begin{equation}\label{TVblow_up2}
		TV(v^{T_n}(\cdot,t),[-\la T_n,\la(T_n+1)])=\f\mbox{ for all }t\in[0,T_n].
		\end{equation}
	\end{enumerate} 
	We use the sequence of data $\{v^{T_n}_0\}$ to define a planar data in multi-D which can be constructed as follows: Let $\{\xi,\zeta_1,\cdots,\zeta_{d-1}\}$ be an orthonormal basis of $\R^d$. Now we define $V_{\xi,T_n}^0$ as follows
	\begin{equation}\label{def:multiD_data}
	V_{\xi,T_n}^0(x):=v_0^{T_n}(x\cdot\xi)\mbox{ for all }x\in\R^d.
	\end{equation} 
	In other words, if $x=(x_1,\cdots,x_d)$ where $x_1=x\cdot\xi,x_{j}=x\cdot\zeta_{j-1}$ for $2\leq j\leq d$, then \eqref{def:multiD_data} can be seen as a function which is constant with respect to $x_i$ variable for $i\geq2$ and it only varies in the first variable, that is $x_1$. It is clear that the function is no more a compact support data. To proceed in a similar way as we have done in Theorem \ref{theorem:1d}, we need a compactly supported initial data. For this purpose let us consider the cube $\mathcal{Q}_r$ in $\R^d$ with each side length $2r$ and centred at origin, that is $\mathcal{Q}_r=[-r,r]^d$ for $r>0$. Define $U_{\xi,T_n}^0$ as follows
	\begin{equation}\label{def:multD_data_cpt}
	U_{\xi,T_n}^0(x)=V_{\xi,T_n}^0(x)\chi_{\mathcal{Q}_{r_n}}(x)\mbox{ where }\chi_{\mathcal{Q}_{r_n}}:=\left\{\begin{array}{ll}
	1&\mbox{ if }x\in\mathcal{Q}_{r_n},\\
	0&\mbox{ otherwise}
	\end{array}\right.
	\end{equation}
	and $r_n=2(\la+M)(T_n+1)$. Let $U_{\xi,T_n}$ be the entropy solution to \eqref{eqn:conlaw} for initial data $U_{\xi,T_n}^0$ as in \eqref{def:multD_data_cpt}. Recall that $v^{T_n}$ is the entropy solution to \eqref{eqn:conlaw} corresponding to initial data $v^{T_n}_0$. Define
	\begin{equation}\label{def:V}
	V_{\xi,T_n}(x,t)=v^{T_n}(x\cdot\xi,t)\mbox{ for }x\in\R^d.
	\end{equation}
	Combining Lemma \ref{lemma:planar_multiD} and Lemma \ref{lemma:linear_transform}, $V_{\xi,T_n}$ is the entropy solution to \eqref{eqn:conlaw} for initial data $V_{\xi,T_n}^0$.  By Proposition \ref{Prop:Kruzkov},
	\begin{align}
	&\int\limits_{|x|\leq \la(T_n+1)+MT_n}|V_{\xi,T_n}(x,t)-U_{\xi,T_n}(x,t)|dx\nonumber\\
	& \leq\int\limits_{|x|\leq \la(T_n+1)+M(T_n+t)}|V_{\xi,T_n}^0(x)-U_{\xi,T_n}^0(x)|dx\label{estimate:contraction}
	\end{align}
	for $t\in[0,T_n]$. Since $r_n=2(\la+M)(T_n+1)$, $R_n:=\la(T_n+1)+2MT_n<r_n$. Subsequently, we get $B(\textbf{0},R_n)\subset \mathcal{Q}_{r_n}$. By \eqref{def:multD_data_cpt},
	\begin{equation*}
	\int\limits_{|x|\leq \la(T_n+1)+M(T_n+t)}|V_{\xi,T_n}^0(x)-U_{\xi,T_n}^0(x)|dx\leq\int\limits_{|x|\leq R_n}|V_{\xi,T_n}^0(x)-U_{\xi,T_n}^0(x)|dx=0
	\end{equation*}
	$\mbox{ for all }t\in[0,T_n]$. Therefore, from \eqref{estimate:contraction} it yields $
	V_{\xi,T_n}(x,t)=U_{\xi,T_n}(x,t)$ for all $(x,t)\in B(\textbf{0},\la(T_n+1)+MT_n)\times[0,T_n]$. From \eqref{def:V} we have
	$ 
	U_{\xi,T_n}(x,t)=v^{T_n}(x\cdot\xi,t)\mbox{ for }(x,t)\in B(\textbf{0},\la(T_n+1)+MT_n)\times[0,T_n].
	$ 
	By Theorem \ref{theorem:BV} and \eqref{TVblow_up2},
	\begin{equation}\label{TV_blowup:Uxi}
	TV(U_{\xi,T_n}(\cdot,t),B(\textbf{0},\la(T_n+1)+MT_n))=\f\mbox{ for }t\in[0,T_n].
	\end{equation}	
	By arranging the sequence of data $\{U_{\xi,T_n}\}$ we want to build a data $U_0$ such that the associated entropy solution $U(\cdot,t)$ does not belong to $BV_{loc}(\R^d)$ for all time $t>0$. For this purpose we use a similar procedure like Theorem \ref{theorem:1d} and arrange them such a way that $U_{\xi,T_n}$ does not interact with $U_{\xi,T_m}$ if $m\neq n$. To this end, define $\rho_n$ as follows	
	$
	\rho_n=\sum\limits_{k=1}^{n}4(T_k+1)(\la+M)+2(T_{n+1}+1)(\la+M).
	$ 
	Now consider initial data $U_0$ and $W_0^m$ defined as follows
	\begin{equation}\label{def:initial_multiD_all_t}
	U_0(x)=\sum\limits_{n=1}^{\f}U_{\xi,T_n}^{0}(x-\rho_n\xi)\mbox{ and }W_0^m(x)=\sum\limits_{n=1}^{m}U_{\xi,T_n}^{0}(x-\rho_n\xi).
	\end{equation}
	Let $U$ and $W^m$ be the entropy solutions to \eqref{eqn:conlaw} corresponding to $U_0$ and $W_0^m$ respectively. We want to show $U(\cdot,t)\notin BV_{loc}(\R^d)$ for all time $t>0$. To do this we prove that $U(\cdot,t)$ agrees with $U_{\xi,T_n}(\cdot-\rho_n\xi,t)$ on some ball around the point $\rho\xi\in\R^d$ for $t\in[0,T_{n}]$ for all $n\geq1$. For this purpose let us first show that $supp(U_{\xi,T_{m+1}}(\cdot-\rho_{m+1}\xi,t))$ is disjoint from $supp(W^m(\cdot,t))$ and $supp (U_{\xi,T_{m+k}}(\cdot-\rho_{m+k}\xi,t))$ for $t\in[0,T_{m+1}]$ and $k\geq2$. 
	By a similar argument as in \eqref{1D:cal1}, from Lemma \ref{lemma:linear_transform},
	\begin{equation}
	\rho_{n-1}+\la(T_{n-1}+1)\leq-\la (T_n+1)+\rho_n.\label{multiD:cal1}
	\end{equation}
	Since $supp(U_{\xi,T_n}^0)$ lies in $B(\textbf{0},\la(T_n+1))$, 
	\begin{equation}\label{support:linear_transform_multiD}
	supp(U_{\xi,T_n}^0(\cdot-\rho_n\xi))\subset B(\rho_n\xi,\la(T_n+1)).
	\end{equation}
	As $\abs{\xi}=1$, from \eqref{def:initial_multiD_all_t}, \eqref{multiD:cal1} and \eqref{support:linear_transform_multiD} we obtain
	\begin{equation*}
	supp(W^m_0)\subset \bigcup\limits_{n=1}^{m}B(\rho_n\xi,\la(T_n+1)) \subset B(\textbf{0},\la(T_n+1)+\rho_n).
	\end{equation*}
	From \eqref{def:initial_multiD_all_t}, we get $\|U_0\|_{L^{\f}(\R)},\|W^m_0\|_{L^{\f}(\R)}\leq u^*$ for all $m\geq1$ due to the fact that $\|U_{\xi,T_n}^0\|_{L^{\f}(\R)}\leq u^*$ and $U^0_{\xi,T_n}$'s have disjoint support for all $n\geq1$. Therefore, $W^m$, the entropy solution to \eqref{eqn:conlaw} with data $W_0^m$ satisfies
	\begin{equation}\label{support:soln_W}
	supp(W^m(\cdot,t))\subset B(\textbf{0},\la(T_n+1)+\rho_n+MT_{m+1})\mbox{ for all }t\in[0,T_{m+1}].
	\end{equation}
	Recall that $U_{\xi,T_{m+1}}$ is the entropy solution to \eqref{eqn:conlaw} with initial data $U^0_{\xi,T_{m+1}}$. By Lemma \ref{lemma:linear_transform}, $U_{\xi,T_{m+1}}(\cdot-\rho_{m+1}\xi,t)$ is the entropy solution to \eqref{eqn:conlaw} associated with initial data $U_{\xi,T_{m+1}}^0(\cdot-\rho_{m+1})$. As we observed before $U_{\xi,T_{m+1}}^0(\cdot-\rho_{m+1}\xi)$ has support in $B(\rho_{m+1}\xi,\la(T_{m+1}+1))$. Hence,
	\begin{equation}\label{support:UT}
	supp(U_{\xi,T_{m+1}}(\cdot-\rho_{m+1}\xi,t))\subset B(\rho_{m+1}\xi,\la(T_{m+1}+1)+MT_{m+1})
	\end{equation}
	for all $t\in[0,T_{m+1}]$. By a similar argument as in \eqref{inequality:cal1},  
	\begin{equation}\label{inequality:cal2}
	-\la( T_{m+1}+1)+\rho_{m+1}-MT_{m+1}>\la(T_m+1)+\rho_m+MT_{m+1}.
	\end{equation}
	Next we wish to show that $B(\rho_{m+1}\xi,\la(T_{m+1}+1)+MT_{m+1})$ and $B(\textbf{0},\la(T_n+1)+\rho_n+MT_{m+1})$ are disjoint. Let $y\in B(\rho_{m+1}\xi,\la(T_{m+1}+1)+MT_{m+1})$. Then $\abs{y-\rho_{m+1}\xi}\leq \la(T_{m+1}+1)+MT_{m+1}$. From the previous observation, $
	\abs{y}\geq \abs{\rho_{m+1}\xi}-\abs{y-\rho_{m+1}\xi}
	\geq\rho_{m+1}-\la(T_{m+1}+1)-MT_{m+1}$.
	Using \eqref{inequality:cal2} we get $
	\abs{y}>\la(T_m+1)+\rho_m+MT_{m+1}$. Hence $B(\rho_{m+1}\xi,\la(T_{m+1}+1)+MT_{m+1})$ and $B(\textbf{0},\la(T_n+1)+\rho_n+MT_{m+1})$ are disjoint. From \eqref{support:soln_W} and \eqref{support:UT},
	\begin{equation}\label{disjoint_support_pre:multi}
	supp(W^m(\cdot,t))\cap supp(U_{\xi,T_{m+1}}(\cdot-\rho_{m+1}\xi,t))=\emptyset\mbox{ for all }t\in[0,T_{m+1}].
	\end{equation}
	Considering $m+k$ instead of $m$ in \eqref{disjoint_support_pre:multi} we conclude that $supp(W^{m+k}(\cdot,t))$ and $supp(U_{\xi,T_{m+k+1}}(\cdot-\rho_{m+k+1}\xi,t))$ are disjoint for $t\in[0,T_{m+k+1}]$. Since $\{T_n\}$ is increasing sequence, $T_{m+k+1}\geq T_{m+1}$ for all $k\geq1$. By previous observation, $supp(U_{\xi,T_{m+k+1}}(\cdot-\rho_{m+k}\xi,t))$ and $supp(W^{m+k}(\cdot,t))$ are disjoint for $t\in[0,T_{m+1}]$ and $k\geq1$. We also know that $supp(U_{\xi,T_{m+1}}(\cdot-\rho_{m+1}\xi,t))$ is contained in $supp(W^{m+k}(\cdot,t))$ for all $k\geq1$. Hence
	\begin{equation}\label{disjoint_support_post:multi}
	supp(U_{\xi,T_{m+k}}(\cdot-\rho_{m+k}\xi,t))\cap supp(U_{\xi,T_{m+1}}(\cdot-\rho_{m+1}\xi,t))=\emptyset
	\end{equation}
	$\mbox{ for all }t\in[0,T_{m+1}]\mbox{ and }k\geq2$. From \eqref{disjoint_support_pre:multi} and \eqref{disjoint_support_post:multi} we conclude that $U_{\xi,T_{m+1}}(\cdot-\rho_{m+1}\xi,t)$ does not interact with $W^{m}$ and $U_{\xi,T_{m+k}}(\cdot-\rho_{m+k}\xi,t)$ for $t\in[0,T_{m+1}]$ and $k\geq2$. Recall that $U$ is the entropy solution to \eqref{eqn:conlaw} with initial data $U$ where $U_0$ is defined as in \eqref{def:initial_multiD_all_t}. Therefore,
	$
	U(x,t)=U_{\xi,T_{m+1}}(x-\rho_{m+1}\xi,t)\mbox{ for }(x,t)\in B(\rho_{m+1}\xi,\la(T_{m+1}+1)+MT_{m+1})\times[0,T_{m+1}].
	$ 
	Recall \eqref{TV_blowup:Uxi} to conclude
		\[
		TV(U(\cdot,t), B(\rho_{m+1}\xi, \la(T_{m+1}+1)+MT_{m+1}))=\f\mbox{ for all }t\in[0,T_{m+1}]\mbox{ and }m\geq1.
		\]
		This completes proof of Theorem \ref{theorem1}.
\end{proof}


\section{Non-existence of $W^{s,p}$-regularizing}\label{sec:Wsp}
With the method of section \ref{sec:TV}, we want to investigate the situation in terms of regularizing in fractional Sobolev space. In this section our main goal is to prove Theorem \ref{theorem2} that is to give an upper bound of $s$ for $L^\f$-$W^{s,p}$ regularizing of entropy solutions to \eqref{eqn:conlaw}. Likewise section \ref{sec:TV} we first show the result for one dimension and then extend it for multi dimension in a similar manner.

\subsection{One dimensional results}\label{sec:Wsp1D}

As we observe in section \ref{sec:TV} that to build a desired entropy solution first and fore most we need a real number sequence $\{a_k\}$ with some properties. Recall that in order to build a non-$BV_{loc}$ entropy solution we were looking for a sequence satisfying \eqref{eq:sequence0} near a point where $f^{\p\p}$ vanishes. Since in this section we are mostly interested in non-degenerate fluxes, we wish to measure the flatness of a flux $f$ near zeros of $f^{\p\p}$ by looking at the ratio $\abs{f^{\p}(a)-f^{\p}(b)}/{\abs{a-b}}$. In purpose of that let us first observe the following.
\begin{lemma}\label{lemma1a}
	Let $p\geq1$. Let $g:[a,b]\rr\re$ be a $C^1$ function for $a<b\in\re$. Fix $y_0\in[a,b]$. Then \textbf{exactly one} of the following occurs:
	\begin{enumerate}
		\item  there is a constant $C>0$ and a sequence $\{y_m\}_{m\geq1}\subset[a,b]$ with $y_m\rr y_0$ such that
		\begin{equation}\label{lm1}
		|g(y_m)-g(y_0)|\leq C|y_m-y_0|^p\mbox{ for all }m\geq1.
		\end{equation}
		\item  for each $\varrho>0$ there exists an $\eta_\varrho>0$ such that
		\begin{equation}\label{eta1}
		\inf\limits_{x\in[a,b],\,0<|x-y_0|<\eta_\varrho}\frac{|g(x)-g(y_0)|}{|x-y_0|^p}\geq \varrho.
		\end{equation}
	\end{enumerate}
\end{lemma}
Note that for fixed $y_0$ and $p$, both of \eqref{lm1}, \eqref{eta1} can not occur simultaneously. Before we see the proof of Lemma \ref{lemma1a}, we want to remark that if $g=f^{\p}$ for some $C^2$ flux $f$, then the exponent $p$ measures the flatness of $f$ near the point $y_0$. Obviously, $p$ and $y_0$ decide which one of \eqref{lm1} and \eqref{eta1} is true. For example consider $g(u)=u^3$. Note that for $p=2$ and $y_0=0$, \eqref{lm1} holds whereas for $p=4$ and $y_0=0$ \eqref{eta1} is true. Observe that for $p=4$ and $y_0=0$, possible choice of $\eta_\varrho$ is $\eta_\varrho=2^{-1}\varrho^{-1}$ when $\varrho\geq 2$ and $\eta_\varrho=4^{-1}$ for $\varrho\in(0,2]$. Similarly, it can be checked that for $p=2$ and $y_0=1$, \eqref{eta1} holds with $\eta_\varrho=\varrho^{-1}$ for $\varrho\geq3$ and $\eta_\varrho=3^{-1}$ for $\varrho\in(0,3]$. It is interesting to note that for $p=3$, $y_0=0$ the inequality \eqref{lm1} is satisfied with $C=1$ and any sequence $\{y_n\}$ such that $y_n\rr0$ for $n\rr\f$. It is clear that for $p=3,y_0=0$, we have $|g(x)-g(y_0)|/|x-y_0|^p=1$, hence it fails to satisfy \eqref{eta1} for all $\varrho>1$.
\begin{proof}
	
	Note that there are two possibilities:
	\begin{description}
		\descitem{Possibility(a):}{Possibility(a)}  for each $\varrho>0$ there exists $\eta_\varrho>0$ such that (\ref{eta1}) holds.
		\descitem{Possibility(b):}{Possibility(b)}  \descref{Possibility(a)}{Possibility(a)} does not hold for some $\varrho_0>0$, i.e. for each $\eta>0$ the inequality (\ref{eta1}) fails to hold for $\varrho_0$.
	\end{description}
	Note that if the \descref{Possibility(a)}{Possibility(a)} holds then there is nothing to prove because it gives (\ref{eta1}). In the \descref{Possibility(b)}{Possibility(b)}, for each $\eta>0$ we have some $x_{\eta}\in[a,b]$ with $0<|x_{\eta}-y_0|<\eta$ such that ${|g(x_{\eta})-g(y_0)|}/{|x_{\eta}-y_0|^p}< \varrho_0$. For a sequence $\eta_m$ such that $\eta_m\rr0$ as $m\rr\f$ we denote $x_m:=x_{\eta_m}$ satisfying ${|g(x_{m})-g(y_0)|}<\varrho_0{|x_{m}-y_0|^p}$ for all $m\geq1$. This gives (\ref{lm1}) with $C=\varrho_0$.
\end{proof}
Our interest lies on behaviour of $f$ near a point where $f^{\p\p}$ vanishes. 
In the next lemma, we consider a flux $f$ such that it is satisfying \eqref{alpha_sup} in an interval on $\R$ and second derivative of $f$ vanishes only at endpoints of the interval then we investigate flatness of $f$ near endpoints of the interval. This lemma plays a key role in proving Theorem \ref{theorem2}.
%
%
%
%
\begin{lemma}\label{lemma:1d:1int}
	Let $a<b$ be two real numbers. Let $f\in C^2([a,b],\re)$ such that $f^{\p\p}(a)=f^{\p\p}(b)=0$ and $f^{\p\p}\neq0$ in $(a,b)$. Suppose $f$ satisfies \eqref{alpha_sup} with $\al_{sup}(f,[a,b])<1$. Let $p$ be a real number such that $1>p^{-1}>\al_{sup}(f,[a,b])$. Then there exist sequences $\{x_k\}_{k\geq1}$, $\{y_k\}_{k\geq1}$ and $M_0>0$ such that \textbf{at least} one of the following holds:
	\begin{enumerate}[(a)]
		\item $x_k\rr a$ as $k\rr\f$ with
		\begin{equation}\label{a}
		\abs{f^{\p}(x_k)-f^{\p}(a)}\leq M_0\abs{x_k-a}^p.
		\end{equation}
		\item $y_k\rr b$ as $k\rr\f$ with 
		\begin{equation}\label{b}
		\abs{f^{\p}(y_k)-f^{\p}(b)}\leq M_0\abs{y_k-b}^p.
		\end{equation}
	\end{enumerate}
	
\end{lemma}
\begin{proof} We invoke Lemma \ref{lemma1a} for $g=f^{\p}$ separately with $y_0=a$ and $y_0=b$. Then we have the following four possibilities:
	\begin{description}
		\descitem{(i)}{(I.)} \eqref{lm1} holds for both $y_0=a$ and $y_0=b$.
		\descitem{(ii)}{(II.)} \eqref{lm1} holds for $y_0=a$ and \eqref{eta1} holds for $y_0=b$.
		\descitem{(iii)}{(III.)} \eqref{eta1} holds for $y_0=a$ and \eqref{lm1} holds for $y_0=b$.
		\descitem{(iv)}{(IV.)} \eqref{eta1} holds for both $y_0=a$ and $y_0=b$.
	\end{description}
	Note that if any of the above possibilities \descref{(I.)}{(i)}-\descref{(III.)}{(iii)} holds then we have nothing to prove because each of possibilities \descref{(I.)}{(i)}-\descref{(II.)}{(iii)} ensures at least one of the inequalities \eqref{a} and \eqref{b}. Next we show that if \descref{(IV.)}{(iv)} holds then it gives a contradiction to the definition of $\al_{sup}(f,[a,b])$. Suppose \descref{(IV.)}{(iv)} holds. Then for $\varrho=2$ there exist $\eta_1,\eta_2>0$ such that the following holds 
	\begin{equation}\label{ab}
	\inf\limits_{x\in(a,a+\eta_1]}\frac{|f^{\p}(x)-f^{\p}(a)|}{|x-a|^p}\geq 2\mbox{ and }	\inf\limits_{x\in[b-\eta_2,b)}\frac{|f^{\p}(x)-f^{\p}(b)|}{|x-b|^p}\geq 2.
	\end{equation}
	Now we claim that there exist ${\de}_1,\de_2\in(0,b-a)$ such that
	\begin{eqnarray}
	&&\inf\limits_{x,y\in[a,a+{\de_1}],x\neq y}\frac{|f^{\p}(x)-f^{\p}(y)|}{|x-y|^p}\geq 1\label{aabb1}\\
	\mbox{ and }&&
	\inf\limits_{x,y\in[b-{\de_2},b],x\neq y}\frac{|f^{\p}(x)-f^{\p}(y)|}{|x-y|^p}\geq 1\label{aabb2}.
	\end{eqnarray}
	We prove only \eqref{aabb1}. The other estimate \eqref{aabb2} follows by a similar argument. If \eqref{aabb1} does not hold for any $\de>0$ then for each $n\in\N$ there exist $x_n,y_n\in[a,a+{n}^{-1}]$ such that $x_n\neq y_n$ and
	$\abs{f^{\p}(x_n)-f^{\p}(y_n)}<\abs{x_n-y_n}^p$. We observe that $x_n,y_n\rr a$ as $n\rr\f$. Without loss of generality suppose $y_n\neq a$ for all $n\geq1$ otherwise we can work with a subsequence of either $x_n$ or $y_n$. As $x_n\rr a$, by using continuity of $f^{\p}$ we get $\abs{f^{\p}(a)-f^{\p}(y_n)}\leq \abs{a-y_n}^p$.	Since $y_n\rr a$ as $n\rr\f$, then $\abs{y_n-a}<\eta_1$ for sufficiently large $n$. This contradicts the first estimate in \eqref{ab}. Therefore, we have \eqref{aabb1}.  Next we show that $f$ satisfies \eqref{condition:LPT} with $\al=p^{-1}$. Since $f^{\p\p}\neq0$ in $(a,b)$ we have
	\begin{equation}\label{B}
	B:=\min\{\abs{f^{\p\p}(x)};\,x\in[a+\de_1,b-\de_2]\}>0.
	\end{equation}
	Set $\bar{\de}=\min\{\de_1,\de_2\}$. Fix $\de\in(0,\min\{\bar{\de}^p,B\abs{a-b+2\bar{\de}}\})$ and $\tau\in\re$. Define	$$A(\de,\tau):=\{x\in[a,b];\,\abs{f^{\p}(x)-\tau}<\de\}.$$ Our next goal is to prove that $\mathcal{L}^1(A(\de,\tau))\leq C\de^{\frac{1}{p}}$ for $\tau\in\R$ and some constant $C>0$. To show this we consider two cases based on the position of $\tau$.
	\begin{description}
		\descitem{Case(1):}{lemma:case1} We first consider the case when $\tau\notin f^{\p}([a,b])$. Since $f^{\p\p}$ does not change its sign in $(a,b)$, $f^{\p}$ is monotone on $[a,b]$. Hence, in this case we have $\min\{\abs{f^{\p}(x)-\tau};\,x\in[a,b]\}$ is either $\abs{f^{\p}(a)-\tau}$ or $\abs{f^{\p}(b)-\tau}$. We deal these two situations separately, 
		\begin{description}
			\descitem{Subcase(1a):}{D1} Suppose $\min\{\abs{f^{\p}(x)-\tau};\,x\in[a,b]\}=\abs{f^{\p}(a)-\tau}$. Then we wish to show
			\begin{equation}\label{setA1}
			A(\de,\tau)\subset\{x\in[a,b];\,\abs{f^{\p}(x)-f^{\p}(a)}<\de\}\subset\{x;\,\abs{x-a}^p<\de\}.
			\end{equation}
			To this end we first observe the following inequality from \eqref{aabb1}
			\begin{equation}\label{a:lw}
			\abs{f^{\p}(x)-f^{\p}(a)}\geq\abs{x-a}^p\mbox{ when }x\in[a,a+\bar{\de}].
			\end{equation}
			Since $f^{\p}$ is monotone in $[a,b]$, we can write
			$$
			\abs{f^{\p}(x)-f^{\p}(a)}=\abs{f^{\p}(x)-f^{\p}(a+\bar{\de})}+\abs{f^{\p}(a+\bar{\de})-f^{\p}(a)}$$ when   $x\in[a+\bar{\de},b-\bar{\de}].
			$
			From \eqref{a:lw} we have 
			\begin{equation}\label{estimate1}
			\abs{f^{\p}(a+\bar{\de})-f^{\p}(a)}\geq\bar{\de}^{p}.
			\end{equation} By using \eqref{B} and Mean Value Theorem for some $c\in(a+\bar{\de},x)$ we have 
			\begin{equation}
			\abs{f^{\p}(x)-f^{\p}(a+\bar{\de})}=\abs{x-a-\bar{\de}}\abs{f^{\p\p}(c)}\geq B\abs{x-a-\bar{\de}}.\label{estimate2}
			\end{equation}
			Combining \eqref{estimate1} and \eqref{estimate2} we have
			\begin{equation*}
			\abs{f^{\p}(x)-f^{\p}(a)}\geq\left\{\begin{array}{lll}
			\abs{x-a}^p&\mbox{ if }&x\in[a,a+\bar{\de}],\\
			\abs{\bar{\de}}^p+B\abs{a+\bar{\de}-x}&\mbox{ if }&x\in[a+\bar{\de},b-\bar{\de}].
			\end{array}\right.
			\end{equation*}
			Note that $f^{\p}$ is monotone. Since $\abs{f^{\p}(a)-f^{\p}(x)}\geq \abs{\de}^p$ holds for $x\in[a+\bar{\de},b-\bar{\de}]$ we have 
			$
			\abs{f^{\p}(a)-f^{\p}(x)}\geq\bar{\de}^p\mbox{ for }x\in[a+\bar{\de},b].
			$ Recall that by our choice of $\de$, we have $\de<\bar{\de}^p$. Therefore, if $\abs{f^{\p}(a)-f^{\p}(x)}<\de$ then we obtain $x\in[a,a+\bar{\de}]$. This shows the second inclusion in \eqref{setA1}. To show the first inclusion of \eqref{setA1} we do the following: fix an $x\in(a,b]$. Suppose $f^{\p}(x)>f^{\p}(a)$. Since $f^{\p}$ is monotone we have $f^{\p}(b)>f^{\p}(a)$. By our assumption we have $\min\{\abs{f^{\p}(x)-\tau};\,x\in[a,b]\}=\abs{f^{\p}(a)-\tau}$ and $\tau\notin[f^{\p}(a),f^{\p}(b)]$. This implies that $f^{\p}(a)>\tau$. Hence we obtain
			\begin{equation}\label{inequality:tau-a}
			\abs{f^{\p}(a)-f^{\p}(x)}\leq\abs{f^{\p}(x)-\tau}.
			\end{equation}
			If $f^{\p}(a)>f^{\p}(x)$ then by a similar argument we can show that \eqref{inequality:tau-a} holds. 
			
			Now if $y\in A(\de,\tau)$ Then we have $\abs{f^{\p}(y)-\tau}<\de$. From \eqref{inequality:tau-a} we have $\abs{f^{\p}(y)-f^{\p}(a)}<\de$. This shows the first inclusion in \eqref{setA1}. Hence \eqref{setA1} is proved.
			\descitem{Subcase(1b):}{D2} Now we consider the case when $\min\{\abs{f^{\p}(x)-\tau};\,x\in[a,b]\}=\abs{f^{\p}(b)-\tau}$. We wish to show
			\begin{equation}\label{setA2}
			A(\de,\tau)\subset\{x\in[a,b];\,\abs{f^{\p}(x)-f^{\p}(b)}<\de\}\subset\{x;\,\abs{x-b}^p<\de\}.
			\end{equation}
			This can be proved in a similar way as we have done for \eqref{setA1} in \descref{D1}{Subcase(1a)}.
		\end{description}
%
		Combining \eqref{setA1} and \eqref{setA2} we get
		\begin{equation}\label{meas(C2)}
		\mathcal{L}^1(A(\de,\tau))\leq 2\de^{\frac{1}{p}}.
		\end{equation}
		\descitem{Case(2):}{lemma:case2} Now we consider the case when $\tau\in f^{\p}([a,b])$. Since $f^{\p\p}$ does not change its sign in $(a,b)$, $f^{\p}$ is monotone on $[a,b]$. Therefore, in this case we have unique $x_0\in[a,b]$ such that $\tau=f^{\p}(x_0)$. By a similar argument as in \descref{lemma:case1}{Case(1)} we observe the following:
		\noi\begin{description}
			\descitem{Subcase(2a):}{case(a)} If $x_0\in[a,a+\bar{\de}]$ then we have
			\begin{equation*}
				\abs{f^{\p}(x)-f^{\p}(x_0)}
				\geq\left\{\begin{array}{lll}
					\abs{x-x_0}^p&\mbox{if }x\in[a,a+\bar{\de}],\\
					\abs{a+\bar{\de}-x_0}^p+B\abs{a+\bar{\de}-x}&\mbox{if }x\in[a+\bar{\de},b-\bar{\de}],\\
					\left(\abs{a+\bar{\de}-x_0}^p+B\abs{a-b+2\bar{\de}}\right.&\mbox{if }x\in[b-\bar{\de},b].\\
					\left.+\abs{b-\bar{\de}-x}^p\right)&
				\end{array}\right.
			\end{equation*}
			Subsequently, by choice of $\de$ we have 
			\begin{equation*}
				A(\de,\tau)	\subset\left\{\begin{array}{lll}
					(x_0-\de^{\frac{1}{p}},x_0+\de^{\frac{1}{p}})&\mbox{if }\de<\abs{x_0-a-\bar{\de}}^p,\\
					(x_0-\de^{\frac{1}{p}},a+\bar{\de}+\de B^{-1})&\mbox{if }\de\geq\abs{x_0-a-\bar{\de}}^p.
				\end{array}\right.
			\end{equation*}
			This yields
			\begin{eqnarray}
			\mathcal{L}^1(A(\de,\tau))	&\leq&\left\{\begin{array}{lll}\label{meas(a)}
			2\de^{\frac{1}{p}}&\mbox{if}&\de<\abs{x_0-a-\bar{\de}}^p,\\
			2\de^{\frac{1}{p}}+\de B^{-1}&\mbox{if}&\de\geq\abs{x_0-a-\bar{\de}}^p.
			\end{array}\right.
			\end{eqnarray}
			\descitem{Subcase(2b):}{case(b)} If $x_0\in[b-\bar{\de},b]$ then we have
			\begin{equation*}
			\abs{f^{\p}(x)-f^{\p}(x_0)}
			\geq\left\{\begin{array}{lll}
					\abs{x-x_0}^p&\mbox{if }x\in[b-\bar{\de},b],\\
					\abs{b-\bar{\de}-x_0}^p+B\abs{b-\bar{\de}-x}&\mbox{if }x\in[a+\bar{\de},b-\bar{\de}],\\
					\left(\abs{a+\bar{\de}-x}^p+B\abs{a-b+\bar{\de}}\right.&\mbox{if }x\in[b-\bar{\de},b].\\
					\left.+\abs{b-\bar{\de}-x_0}^p\right)&
				\end{array}\right.
			\end{equation*}
			By a similar argument as in \descref{case(a)}{Subcase(2a)} we have
			\begin{eqnarray}
			\mathcal{L}^1(A(\de,\tau))	&\leq&\left\{\begin{array}{lll}\label{meas(b)}
			2\de^{\frac{1}{p}}&\mbox{if}&\de<\abs{b-\bar{\de}-x_0}^p,\\
			2\de^{\frac{1}{p}}+\de B^{-1}&\mbox{if}&\de\geq\abs{b-\bar{\de}-x_0}^p.
			\end{array}\right.
			\end{eqnarray}
			\descitem{Subcase(2c):}{case(c)} If $x_0\in[a+\bar{\de},b-\bar{\de}]$ we have
			\begin{equation*}
			\abs{f^{\p}(x)-f^{\p}(x_0)}		\geq\left\{\begin{array}{lll}
			B\abs{x_0-x}&\mbox{if }x\in[a+\bar{\de},b-\bar{\de}],\\
			B\abs{a+\bar{\de}-x_0}+\abs{a+\bar{\de}-x}^p&\mbox{if }x\in[a,a+\bar{\de}],\\
			B\abs{b-\bar{\de}-x_0}+\abs{b-\bar{\de}-x}^p&\mbox{if }x\in[b-\bar{\de},b].
			\end{array}\right.
			\end{equation*}
			Subsequently, we have
			\begin{equation*}
				A(\de,\tau)\subset\left\{\begin{array}{lll}
					\left(x_0-\frac{\de}{ B},x_0+\frac{\de}{ B}\right)&\mbox{if }\abs{a+\bar{\de}-x_0},\abs{b-\bar{\de}-x_0}\geq\frac{\de}{B},\\
					\left(a+\bar{\de}-\de^\frac{1}{p},x_0+\frac{\de}{ B}\right)&\mbox{if }\abs{a+\bar{\de}-x_0}<\frac{\de}{ B}\leq\abs{b-\bar{\de}-x_0},\\
					\left(x_0-\frac{\de}{ B},b-\bar{\de}+\de^{\frac{1}{p}}\right)&\mbox{if }\abs{b-\bar{\de}-x_0}<\frac{\de}{ B}\leq\abs{a+\bar{\de}-x_0},\\
					\left(a+\bar{\de}-\de^\frac{1}{p},b-\bar{\de}+\de^{\frac{1}{p}}\right)&\mbox{if }\abs{a+\bar{\de}-x_0},\abs{b-\bar{\de}-x_0}<\frac{\de}{ B}.
				\end{array}\right.
			\end{equation*}
			This yields 
			\begin{equation}
			\mathcal{L}^1(A(\de,\tau))\leq\left\{\begin{array}{lll}\label{meas(c)}
			{2\de}{B}^{-1}&\mbox{if }\abs{a+\bar{\de}-x_0},\abs{b-\bar{\de}-x_0}\geq\frac{\de}{B},\\
			\de^\frac{1}{p}+{2\de}{B}^{-1}&\mbox{if }\abs{a+\bar{\de}-x_0}<\frac{\de}{ B}\leq\abs{b-\bar{\de}-x_0},\\
			\dip	\de^\frac{1}{p}+{2\de}{B}^{-1}&\mbox{if }\abs{b-\bar{\de}-x_0}<\frac{\de}{ B}\leq\abs{a+\bar{\de}-x_0},\\
			2\de^\frac{1}{p}+{2\de}{B}^{-1}&\mbox{if }\abs{a+\bar{\de}-x_0},\abs{b-\bar{\de}-x_0}<\frac{\de}{ B}.
			\end{array}\right.
			\end{equation}
		\end{description}
		Combining \eqref{meas(a)}, \eqref{meas(b)} and \eqref{meas(c)} we have
		\begin{equation}\label{meas(C1)}
		\mathcal{L}^1(A(\de,\tau))\leq 2\de^{\frac{1}{p}}+2B^{-1}\de.
		\end{equation}
	\end{description}
	From \eqref{meas(C1)} and \eqref{meas(C2)} we have $\mathcal{L}^1(A(\de,\tau))\leq 2\de^{\frac{1}{p}}+2B^{-1}\de\leq\left(2+2B^{-1}\right)\de^{\frac{1}{p}}$. Define $$C_1:=\max\left\{(b-a)(\min\{\bar{\de}^p,B\abs{a-b+2\bar{\de}}\})^{-\frac{1}{p}},\left(2+2B^{-1}\right)\right\}.$$ Note that $C_1$ does not depend on $\de$ and $\tau$. Now observe that for $\de>\min\{\bar{\de}^p,B\abs{a-b+2\bar{\de}}\}$ and $\tau\in\re$ we have
	$\mathcal{L}^1(\{x\in[a,b];,\abs{f^{\p}(x)-\tau}<\de\})\leq b-a\leq C_1\de^{\frac{1}{p}}$. Therefore, we have $\mathcal{L}^1(A(\de,\tau))\leq C_1\de^{\frac{1}{p}}\mbox{ for all }\de>0,\,\tau\in\re$. This shows that $f$ satisfies \eqref{condition:LPT} for $\al=1/{p}$. This contradicts the definition of $\al_{sup}(f,[a,b])$ since ${1}/{p}>\al_{sup}(f,[a,b])$.
\end{proof}

Next we want to remark few things about fluxes satisfying \eqref{condition:LPT} for some $\al\in(0,1]$ and $C>0$. First we recall the following fact implicitly mentioned in \cite{Lions1}.
\begin{remark}\label{remark:nonzero}
	If $f\in C^2(I,\R)$ be a 1-D flux satisfying \eqref{condition:LPT} for some $\al\in(0,1]$ then $f^{\p\p}$ can not vanish in a non-trivial sub-interval of $I$. Note that if $f^{\p\p}$ vanishes on some $I_{1}\subset I$ then $f^{\p}$ is constant on $I_1$. Hence if we choose $\tau=f^{\p}(a_0)$ for some $a_0\in I_1$ then we get $I_1\subset\{a\in I;\abs{f^{\p}(a)-\tau}<\de\}\mbox{ for any }\de>0$.	This violates \eqref{condition:LPT}. 
\end{remark}
To give more insight on $\al_{sup}$ we make a simple observation on non-degeneracy of uniformly convex flux in the following remark.
\begin{remark}\label{remark:uniformlyconvex}
	Let $I=[a,b]$ for some $a<b$ and $f\in C^2(I,\R)$. Then $f$ satisfies \eqref{condition:LPT} with $\al=1$ if and only if $f^{\p\p}\neq0$ on $I$.
\end{remark}
\begin{proof}
	Suppose $f$ satisfies \eqref{condition:LPT} with $\al=1$. Suppose there exists a point $a_0\in I$ such that $f^{\p\p}(a_0)=0$. By Remark \ref{remark:nonzero} there $f^{\p\p}$ can not vanish identically on $I$. Therefore, we can find a point $\bar{a}$ such that $f^{\p\p}\neq0$ either on $(\bar{a},\bar{a}+\e)$ or on $(\bar{a}-\e,\bar{a})$ but $f^{\p\p}(\bar{a})=0$. Without loss of generality we assume that $f^{\p\p}\neq0$ on $(\bar{a},\bar{a}+\e)$. Consider a sequence $\{b_n\}\subset (\bar{a},\bar{a}+\e)$ such that $b_n\rr \bar{a}$ as $n\rr\f$. Choose $\tau_n=f^{\p}(b_n)$ and $\de_n=2\abs{f^{\p}(\bar{a})-f^{\p}(b_n)}$. Since $f^{\p\p}\neq0$ on $(\bar{a},\bar{a}+\e)$ we have $f^{\p}$ is monotone on $(\bar{a},\bar{a}+\e)$. Hence we have $[\bar{a},b_n]\subset\{u\in I;\,\abs{f^{\p}(u)-\tau}<\de_n\}$. From \eqref{condition:LPT} we obtain $\frac{f^{\p}(\bar{a})-f^{\p}(b_n)}{\abs{\bar{a}-b_n}}\geq\frac{1}{2C}$. Passing to the limit as $n\rr\f$ we have $f^{\p\p}(\bar{a})\neq0$. This is a contradiction. Hence $f^{\p\p}\neq0$ on $I$. Conversely, suppose $f^{\p\p}\neq0$ on $I$. Let $m=\min\{\abs{f^{\p\p}(u)};\,u\in I\}$. Then we have $\abs{f^{\p}(u)-f^{\p}(v)}\geq m\abs{u-v}$.	By using a similar technique as in Lemma \ref{lemma:1d:1int} we can show that $f$ satisfies \eqref{condition:LPT} with $\al=1$.
	\end{proof}
By Remark \ref{remark:uniformlyconvex} we can conclude that for uniformly convex flux $\al_{sup}=1$. It is crucial to note that converse is not true, i.e. there are fluxes in one dimension which are not uniformly convex but still satisfies \eqref{alpha_sup} with $\al_{sup}=1$. Next we give an example of such flux.
\begin{example}\label{example:logflux}
	Let $f:[0,e^{-1}]\rr\R$ be defined as follows
	$
	f(u):=\int\limits_{0}^{u}\int\limits_{0}^{y}\frac{1}{\abs{\log x}}\,dxdy.
	$
\end{example}
Now we claim:
\begin{lemma}\label{lemma:logflux}
	Let $f$ be defined as in Example \ref{example:logflux}. Then $\al_{sup}(f,[0,e^{-1}])=1$. 
\end{lemma}
\begin{proof}
 It is clear that $f\in C^2([0,e^{-1}])$. To show $\al_{sup}(f,[0,e^{-1}])$, we first prove that $f$ satisfies \eqref{condition:LPT} with $\al=k/(k+1)$ for each $k\geq 2$. From definition of $f$, we observe that $f^{\p}(u)=\int\limits_{0}^{u}\frac{1}{\abs{\log x}}\,dx$.
 Fix a $k\geq2$. From \eqref{inequality:log} we have for $u\in[0,e^{-k^2}]$
 $
 f^{\p}(u)\geq \int\limits_{0}^{u}\abs{x}^{\frac{1}{k}}\,dx=\frac{k}{k+1}u^{1+{k}^{-1}}\geq2^{-1}u^{1+k^{-1}}.
 $
 Since $f^{\p\p}=\abs{\log x}^{-1}\geq0$, we have $f^{\p}$ is increasing and $f^{\p}(0)=0$. Hence we have $\abs{f^{\p}}\geq 2^{-1}u^{1+k^{-1}}.$ Note that for $u\in[e^{-k^2},e^{-1}]$, $\abs{f^{\p\p}}\geq k^{-2}$. By a similar technique as in proof of Lemma \ref{lemma:1d:1int} we show that $f$ satisfies \eqref{condition:LPT} with $\al=k/(k+1)$ and some constant $C$ (depends on $k$). Recall the definition of $\al_{sup}$ as in \eqref{alpha_sup} and observe that $\sup\{k/(k+1);\,k\geq1\}=1$. Therefore, we conclude that the function $f$ as in Example \ref{example:logflux} satisfies $\al_{sup}(f,[0,e^{-1}])=1$.
\end{proof}
In Lemma \ref{lemma:logflux} we have seen that in one dimension $\al_{sup}$ can be 1 for certain $C^2$ fluxes. In the next lemma we show that for a slightly regular flux in multi-D i.e. $F\in C^{2,\B}(I,\R^d)$ the non-degeneracy exponent $\al_{sup}$ is strictly smaller than 1. This plays a key role to prove Proposition \ref{proposition2}.
\begin{lemma}\label{lemma:Holder}
	Let $d>1$ and $F\in C^{2,\B}([a,b],\re^d)$ for some $0<\B\leq1$, $a<b$. Let $\al_{sup}(F,[a,b])$ be determined by (\ref{alpha_sup}). Then $\al_{sup}(F,[a,b])\leq (1+\B)^{-1}$. Moreover, there exist $x_0\in I$, $\xi_0\in \mathcal{S}^{d-1}$, $r>0$, $C_0>0$ and $\{x_k\}_{k\geq1}\subset(x_0,x_0+r)\subset I$ such that $F^{\p\p}\cdot\xi_0\neq0$ on $(x_0,x_0+r)$, $x_k\rr x_0$ as $k\rr\f$ and 
	\begin{equation}\label{sequence:Holder}
	\abs{\xi_0\cdot F^{\p}(x_k)-\xi_0\cdot F^{\p}(x_0)}\leq C_0\abs{x_k-x_0}^{1+\B}\mbox{ for all }k\in\N.
	\end{equation}	
\end{lemma}
\begin{proof}
	Since $F\in C^{2,\B}([a,b])$, there exists a constant $C$ such that $
	|F^{\p\p}(x)-F^{\p\p}(y)|\leq C|x-y|^{\B}$ for all $x,y\in[a,b]$. Fix a $x_0\in[a,b]$. We choose a direction $\xi_0\in\mathcal{S}^{d-1}$ such that $\xi_0\cdot F^{\p\p}(x_0)=0$. Therefore, by Mean Value Theorem we have
	\begin{eqnarray}
	|\xi_0\cdot F^{\p}(x)-\xi_0\cdot F^{\p}(x_0)|&=&|\xi_0\cdot F^{\p\p}(y)(x-x_0)|,\hspace{.05cm}\mbox{ for some $y$ between $x$ and $x_0$}\nonumber\\
	&=&|\xi_0\cdot F^{\p\p}(y)-\xi_0\cdot F^{\p\p}(x_0)||x-x_0|\nonumber\\
	&\leq&C|y-x_0|^{\B}|x-x_0|\nonumber\\
	&\leq&C|x-x_0|^{1+\B}\hspace{0.5cm}\mbox{ since, }|y-x_0|\leq|x-x_0|.\label{fprime}
	\end{eqnarray}
	If $\al_{sup}(F,[a,b])>1/(1+\B)$ then there exists an $\al_0\in(1/(1+\B),\al_{sup})$ such that $F$ satisfies (\ref{condition:LPT}) with $\al_0$ and some $C_1>0$. 
	Consider $\tau=-\xi_0\cdot F^{\p}(x_0)$. Define $\tau_1:=C_0^{-1}\tau\mbox{ and }\xi_1:=C_0^{-1}\xi_0$ with $C_0=\left(\tau^2+\abs{\xi_0}^2\right)^{1/2}$. This yields $\tau_1^2+\abs{\xi_1}^2=1$. Since $F$ satisfies \eqref{condition:LPT} with exponent $\al_0$ we have $
	\mathcal{L}^1\left(\left\{x\in[a,b];|\tau_1+\xi_1\cdot F^{\p}(x)|<{\de}C_0^{-1}\right\}\right)\leq C_1\left({\de}C_0^{-1}\right)^{\al_0}$ for $\de\in(0,1)$. Therefore, we get
	\begin{equation}\label{alpha_0}
	\mathcal{L}^1\left(\left\{x\in[a,b];|\xi_0\cdot F^{\p}(x_0)-\xi_0\cdot F^{\p}(x)|<\de\right\}\right)\leq \frac{C_1}{C_2}\de^{\al_0}
	\end{equation}
	where $C_2=C_0^{\al_0}$. On the other hand estimate (\ref{fprime}) gives
	\begin{equation*}
	\mathcal{L}^1\left(\left\{x\in[a,b];|\tau+\xi_0\cdot F^{\p}(x)|<\de\right\}\right)\geq  \frac{1}{C^{{1}/({1+\B})}}\de^{{1}/({1+\B})}.
	\end{equation*}
	By (\ref{alpha_0}) we have
	\begin{eqnarray}
	\frac{C_1}{C_2}\de^{\al_0}\geq \frac{1}{C^{{1}/({1+\B})}}\de^{{1}/({1+\B})}
	\mbox{ which implies, }\de^{\al_0-1/(1+\B)}\geq\frac{C_2}{C_1C^{{1}/({1+\B})}}.\label{delta1}
	\end{eqnarray}
	Since $\al_0>1/(1+\B)$ and (\ref{delta1}) is true for all $\de\in(0,1)$ this gives a contradiction. Note that \eqref{sequence:Holder} is guaranteed by \eqref{fprime}. This completes the proof of Lemma \ref{lemma:Holder}.
\end{proof}

Now we are almost ready to prove Theorem \ref{theorem2}. Before that we wish to discuss more on the regularity assumption \eqref{mild_regularity}. In one dimension, if $f$ satisfies \eqref{alpha_sup} and $f^{\p\p}$ vanishes finitely many times in an interval $[a,b]$ then it satisfies \eqref{mild_regularity}. We have the following remark for $C^\f$ fluxes.
\begin{remark}\label{lemma:C-infinity}
	Let $I=[a,b]$ for $-\f<a<b<\f$ and $F\in C^\f(I,\R^d)$ be satisfying \eqref{alpha_sup} with $\al_{sup}(F,I)>0$. Then $F$ satisfies \eqref{mild_regularity}.
\end{remark}
\begin{proof}
	First we observe that for $d=1$, if $f\in C^\f(I,\R)$ satisfies \eqref{alpha_sup} with $\al_{sup}>0$ then $f^{\p\p}$ has finitely many zeros in $I$ (for instance see \cite{Ju14}). Since $I$ is compact, if $f^{\p\p}$ has infinitely many zeros then there exists a point $\bar{a}\in I$ such that $f^{(n)}(\bar{a})=0$ for all $n\geq2$. This shows that $f$ can not satisfy \eqref{condition:Ju14} with $d_f<\f$. This violates the assumption $\al_{sup}(f,I)=d_f^{-1}>0$. Hence $f^{\p\p}$ has finitely many zeros in $I$. Therefore, by previous observation, we conclude that $f$ satisfies \eqref{mild_regularity}. Consequently, $u\mapsto F^{\p\p}(u)\cdot\xi$ satisfies \eqref{mild_regularity} for each $\xi\in\mathcal{S}^{d-1}$. Hence to show that $F$ satisfies \eqref{mild_regularity} it is enough to prove the following
	\begin{equation}\label{cal:smooth1}
	\al_{sup}(F,I)=\inf\{\al_{sup}(F\cdot\xi,I),\xi\in\mathcal{S}^{d-1}\}.
	\end{equation}
	Since $F,F\cdot\xi$ are smooth functions, $d_{F}[I],d_{F\cdot\xi}[I]$ are positive integers (for more detail see \cite{Ju14}). By previous observation we know that \eqref{cal:smooth1} with $\al_{sup}>0$ is equivalent to $d_F[I]=\sup\{d_{F\cdot\xi}[I],\,\xi\in\mathcal{S}^{d-1}\}<\f$. Note that $d_{F}[I]\geq d_{F\cdot\xi}[I]$ for each $\xi\in\mathcal{S}^{d-1}$. To show this fix a $k\geq d_F[I]$. From \eqref{condition:Ju14} we have $\mbox{span}\{F^{\p\p}(u),\cdots,F^{(k+1)}(u)\}=\R^d$ for all $u\in I$. For a $\xi\in\mathcal{S}^{d-1}$ we have $\mbox{span}\{F^{\p\p}(u)\cdot\xi,\cdots,F^{(k+1)}(u)\cdot\xi\}=\R$ for all $u\in I.
	$ Hence we have $d_{F\cdot\xi}[I]\leq k$ for each $\xi\in\mathcal{S}^{d-1}$. Taking infimum over all $k\geq d_F[I]$ we have $d_{F\cdot\xi}\leq d_F[I]$. Therefore, we obtain $d_F[I]\geq\sup\{d_{F\cdot\xi}[I],\,\xi\in\mathcal{S}^{d-1}\}$. Now we show that $d_F[I]\leq \sup\{d_{F\cdot\xi}[I],\,\xi\in\mathcal{S}^{d-1}\}$.	Suppose $d_{F\cdot\xi}[I]\leq d_F[I]-1$ for all $\xi\in\mathcal{S}^{d-1}$. Then $\{F^{\p\p}(u)\cdot\xi,\cdots,F^{(d_F[I])}(u)\cdot\xi\}$ spans $\R$ for all $u\in I$ and $\xi\in\mathcal{S}^{d-1}$. Hence $\{F^{\p\p}(u),\cdots,F^{(d_F[I])}(u)\}$ spans whole $\R^d$ for all $u\in I$. This gives a contradiction with definition \eqref{condition:Ju14} of $d_F[I]$.
\end{proof}
In the following example we construct a $C^2$ flux in 1-D such that $f^{\p\p}$ vanishes at infinitely many points in a finite interval and it also satisfies \eqref{mild_regularity}. 
\begin{example}\label{example:infinitely}
	Here we provide a non-trivial example of flux in 1-D satisfying \eqref{mild_regularity} and zero set of second derivative is not locally finite. Define $g_n:[1/(n+1),1/n]\rr\re$ for $n\geq2$, as follows:
	\begin{equation*}
	g_n(x):=\left\{\begin{array}{lll}
	a_n\left(x-(n+1)^{-1}\right)^{1-{n}^{-1}}&\mbox{if}&x\in\left[\frac{1}{n+1},\frac{1}{n+1}+\frac{1}{2(n+1)^{2}}\right],\\
	{n}^{-1}&\mbox{if}&x\in\left[\frac{1}{n+1}+\frac{1}{2(n+1)^{2}},\frac{1}{n}-\frac{1}{2(n+1)^{2}}\right],\\
	a_n\left({n}^{-1}-x\right)^{1-{n}^{-1}}&\mbox{if}&x\in\left[\frac{1}{n}-\frac{1}{2(n+1)^{2}},\frac{1}{n}\right]
	\end{array}\right.
	\end{equation*}
	where $a_n:=2^{1-\frac{1}{n}}\frac{(n+1)^{2-\frac{2}{n}}}{n}$.	Suppose $g:[0,\frac{1}{2}]\rr\re$ defined as $g:=\sum\limits_{n=2}^{\f}\chi_{[\frac{1}{n+1},\frac{1}{n}]}g_n$. Note that $g$ is a continuous function vanishing on $\{1/n\}_{n\geq2}\cup\{0\}$. We define a flux $f:[0,1/2]\rr\re$ as 
	$f(x)=\int\limits_{0}^{x}\int\limits_{0}^{y}g(z)\,dzdy.$
\end{example}
Then we claim:
\begin{lemma}
	Let $f$ be defined as in Example \ref{example:infinitely}. Then the following holds
	\begin{enumerate}
		\item $\al_{sup}\left(f,\left[0,{1}/{2}\right]\right)={1}/{2}$ and for each $n\geq1$ we have $\al_{sup}\left(f,\left[{1}/({n+1}),{1}/{n}\right]\right)=\frac{n}{2n-1}>\frac{1}{2}$.
		\item $f$ satisfies the mild regularity assumption \eqref{mild_regularity}, i.e.,
		\begin{equation}\label{mild_reg:ex}
		\al_{sup}\left(f,\left[0,{1}/{2}\right]\right)=\inf\left\{\al_{sup}\left(f,\left[{1}/{(n+1)},{1}/{n}\right]\right),\,n\geq2\right\}.
		\end{equation}
	\end{enumerate}
\end{lemma}
\begin{proof}Since $g\geq0$ we have $f^{\p}$ is non-decreasing. Observe that
	\begin{equation*}
	\frac{1}{3}\frac{1}{n(n+1)^2}+\frac{1}{n^2(n+1)^2}\leq f^{\p}\left({n}^{-1}\right)-f^{\p}\left((n+1)^{-1}\right)\leq \frac{1}{2}\frac{1}{n(n+1)^2}+\frac{1}{n^2(n+1)^2}.
	\end{equation*}
	Let $k,l\in\N$ such that $k<l$. We consider $x\in(1/(l+1),1/l)$ and $y\in(1/(k+1),1/k)$, then we have
	\begin{eqnarray*}
		f^{\p}(y)-f^{\p}(x)&\geq& \frac{1}{6}\sum\limits_{n=k+1}^{l-1}\left[\frac{2}{n(n+1)^2}+\frac{1}{n^2(n+1)^2}\right]\\
		&\geq&
		\frac{1}{6}\sum\limits_{n=k+1}^{l-1}\left[\frac{1}{n^2}-\frac{1}{(n+1)^2}\right]
		\geq \frac{1}{6}\left[\frac{1}{(k+1)^2}-\frac{1}{l^2}\right]
		\geq \frac{1}{6}\left(\frac{l-k-1}{(k+1)l}\right)^2.
	\end{eqnarray*}
	Observe that if $
	\left[(l+1)^{-1},k^{-1}\right]\subset\{x\in[0,1];\abs{f^{\p}(x)-f^{\p}(x_0)}<\de\}\subset\left[\frac{1}{l+2},\frac{1}{k-1}\right]
	$ for $l-k\geq2$, then we have $
	\frac{l-k-1}{(k+1)l}\leq(12\de)^{\frac{1}{2}}$.	Therefore
	\begin{equation*}
	\frac{1}{k-1}-\frac{1}{l+2}
	=\left(\frac{k+1}{k-1}\right)\left(\frac{l}{l+2}\right)\left(\frac{l+3-k}{l-1-k}\right)\left(\frac{l-k-1}{(k+1)l}\right)	\leq 15(12\de)^{\frac{1}{2}}.
	\end{equation*}
	Fix  $x_0\in\left[\frac{1}{n+1},\frac{1}{n+1}+\frac{1}{2(n+1)^{2}}\right]$. Then observe the following cases:
	\begin{enumerate}
		\item if $x\in\left[\frac{1}{n+1},\frac{1}{n+1}+\frac{1}{2(n+1)^{2}}\right]$ then we have
		$
		f^{\p}(x)-f^{\p}(x_0)\geq \frac{a_n}{3}(x-x_0)^{2-\frac{1}{n}}.
		$
		\item if $x\in\left[(n+1)^{-1}+2^{-1}(n+1)^{-2},{n}^{-1}-2^{-1}(n+1)^{-2}\right]$ then we have
		$$
		f^{\p}(x)-f^{\p}(x_0)\geq\frac{a_n}{3}\left(\frac{1}{n}+\frac{1}{2(n+1)^{2}}-x_0\right)^{2-\frac{1}{n}}+\left(x-\frac{1}{n+1}-\frac{1}{2(n+1)^{2}}\right)^{\frac{4}{3}}.
		$$
		\item if $x\in\left[\frac{1}{n}-\frac{2(n+1)^{2}},\frac{1}{n}\right]$ then we have
		\begin{equation*}
		\begin{array}{lll}
		f^{\p}(x)-f^{\p}(x_0)&\geq&\frac{a_n}{3}\left(\frac{1}{n}+\frac{1}{2(n+1)^{2}}-x_0\right)^{2-\frac{1}{n}}+\left(n(n+1)^2\right)^{-\frac{4}{3}}\\
		&+&\frac{a_n}{3}\left(x-\frac{1}{n}+\frac{1}{2(n+1)^{2}}\right)^{2-\frac{1}{n}}.
		\end{array}
		\end{equation*}
	\end{enumerate}
	By a similar argument as in proof of Lemma \ref{lemma:1d:1int} we can show that
	\begin{equation}\label{meas:ex}
	\mathcal{L}^1\left(\{x\in[0,1/2];\,\abs{f^{\p}(x)-f^{\p}(x_0)}<\de\}\right)\leq C\de^{\frac{1}{2}}.
	\end{equation}
	For other positions of $x_0$ we can show \eqref{meas:ex} by a similar manner. Therefore, it shows that $f$ satisfies \eqref{condition:LPT} with $\al={1}/{2}$ on the interval $[0,1/2]$. Also note that in a neighbourhood of $1/(n+1)$, $f^{\p\p}$ is Holder continuous with exponent $\frac{n-1}{n}$. By Lemma \ref{lemma:Holder} we have $\al_{sup}(f,[0,1])\leq\frac{n}{2n-1}$. Therefore, we have $	\al_{sup}\left(f,\left[0,{1}/{2}\right]\right)=\frac{1}{2}$. By a similar argument as in proof of Lemma \ref{lemma:1d:1int} we can show that $\al_{sup}\left(f,\left[(n+1)^{-1},{n}^{-1}\right]\right)=\frac{n}{2n-1}$. Hence we prove \eqref{mild_reg:ex}

\end{proof}
Now we are ready to proceed for proof of Theorem \ref{theorem2}. Similar to section \ref{sec:TV} here we first construct an entropy solution $u$ such that $u\notin W^{s,1}_{loc}$ up to a finite time $t_0$. This is one of the key tools to prove Theorem \ref{theorem2}.
\begin{proposition}\label{Proposition:1dW_uptoT}
	Let $R_0>0$. Let $f\in C^2([-R_0,R_0],\re)$ be a flux satisfying the non-degeneracy condition (\ref{alpha_sup}) with $\al_{sup}(f,[-R_0,R_0])\in(0,1)$ and additionally \eqref{mild_regularity}. Then there exist constants $\la>0$ and $u^*>0$ such that the following holds: for every $T>0$ and $s>\al_{sup}(f,[-R_0,R_0])$ there exists an initial data $v_0^T$ verifying the following properties: 
	\begin{enumerate}
		\item $\|v^T_0\|_{L^{\f}(\R)}\leq u^*$ and $supp(v^T_0)\subset[-\la T,\la(T+1)]$.
		\item if $v^T$ is the entropy solution to \eqref{eqn:conlaw} with initial data $v^T_0$, then $v^T(\cdot,t)\notin W^{s,1}([-\la T,\la(T+1)])$ for all $t\in[0,T]$.
	\end{enumerate}
\end{proposition}
\begin{proof}
	Fix an $s\in(\al_{sup}(f,[-R_0,R_0]),1)$. Now we choose a $p_0$ such that \\$\al_{sup}(f, [-R_0,R_0])<p_0^{-1}<s$. Since $\al_{sup}(f,[-R_0,R_0])<1$ the zero set of $f^{\p\p}$, $\mathcal{Z}(f,[-R_0,R_0])$ is non-empty due to Remark \ref{remark:uniformlyconvex}. As $\al_{sup}>0$ $\mathcal{Z}(f,[-R_0,R_0])$ has empty interior by Remark \ref{remark:nonzero}. Therefore, $[-R_0,R_0]\setminus\mathcal{Z}(f,[-R_0,R_0])$ is non-empty open set and it can be written as at most countable union of disjoint open intervals, that is,
	$
	[-R_0,R_0]\setminus\mathcal{Z}(f,[-R_0,R_0])=\bigcup\limits_{j\in \mathcal{J}}I_j
$
	for some index set $\mathcal{J}\subset\mathbb{N}$. Since $f$ satisfies \eqref{mild_regularity}, we have
	$
	\al_{sup}(f,[-R_0,R_0])=\inf\limits_{j\in\mathcal{J}}\al_{sup}(f,I_j).
	$ Hence there exists an interval $I_{j_0}=[a_0,b_0]\subset\R$ such that $\al_{sup}(f,[a_0,b_0])<p_0^{-1}$. Since $I_j$ are disjoint interval and subset of $[-R_0,R_0]\setminus\mathcal{Z}(f,[-R_0,R_0])$ we have $f^{\p\p}\neq0$ in $(a_0,b_0)$ and $f^{\p\p}(a_0)=f^{\p\p}(b_0)=0$. Now we invoke Lemma \ref{lemma:1d:1int} for interval $[a_0,b_0]$ and show that at least one of \eqref{a} and \eqref{b} holds. Without loss of generality we assume that \eqref{a} holds. Hence there exists a sequence $\{a_k\}_{k\geq1}$ such that 
	\begin{equation}\label{theorem2:asmp1}
	|f^{\p}(a_k)-f^{\p}(a_0)|\leq C|a_k-a_0|^{p_0}\mbox{ for all }k\geq1.
	\end{equation}
	Without loss of generality we assume that $C=1$. Fix a time $T>0$. We wish to construct an initial data in a similar way as we have done in Proposition \ref{Proposition:1d_uptoT}. Note that \descref{linear}{Case(1)} of Proposition \ref{Proposition:1d_uptoT} does not appear here because of Remark \ref{remark:nonzero}. Therefore, to construct a data we only need to mimic the construction in \descref{non-linear}{Case(2)} of Proposition \ref{Proposition:1d_uptoT}. To set the ground we first invoke Lemma \ref{lemma:subsequence}, to get a subsequence of $\{a_k\}$ (still labelled as $\{a_k\}$) verifying the following properties:
	\begin{equation}\label{WLOG:akqk}
	|a_k-a_0|=k^{-q_{k}}\mbox{ for some sequence }q_k\mbox{ satisfying }q_{k+1}\geq q_k+1\mbox{ for all }k\geq1.
	\end{equation}
	We define $\e_k,A_k$ as follows 
	\begin{equation}\label{theorem2:asmp2}
	\e_k=|a_k-a_0|^{p_0}\mbox{ and }A_k=\frac{1}{k^{1+\e}}\mbox{ for some }\e\in(0,1).
	\end{equation}
	Since $q_k\geq1$ and $p_0>1$ we can choose $\e>0$ small enough such that $p_0q_k-1-\e>0$. By using \eqref{WLOG:akqk} and \eqref{theorem2:asmp2} we get $
	A_k\e_k^{-1}=k^{-1-\e}\abs{a_k-a_0}^{-1}=k^{-1-\e}k^{q_kp_0}\geq1$ since $p_0q_k-1-\e>0$.
	By employing \eqref{theorem2:asmp1}, \eqref{theorem2:asmp2} we have $\abs{f^{\p}(a_k)-f^{\p}(a_0)}\leq \e_k\leq A_k$ and $\sum\limits_{k\geq1}A_k+2\e_1<\f$. This ensures \eqref{condition:ekAk}. Recall the definitions of $N_k,J_k,\si_m,\de_n,x_n$ as in \eqref{def:NkJk}, \eqref{def:sigma}, \eqref{def:delta} and \eqref{def:x_n} respectively. Consider the initial data $w_0$ defined as in \eqref{def:initial1} with $\e_k,A_k$ are as in \eqref{theorem2:asmp2}. Similar to Proposition \ref{Proposition:1d_uptoT}, there exists a constant $u^*$ such that $	\|w_0\|_{L^{\f}(\R)}\leq u^{*}$.  Let $M$ be defined as 
	\begin{equation}\label{def:Wsp:M}
	M:=\max\{\abs{f^{\p}(u)};\,\abs{u}\leq u^*\}.
	\end{equation}
    Let $w$ be the entropy solution to \eqref{eqn:conlaw} with initial data $w_0$. As we have seen in proof of Proposition \ref{Proposition:1d_uptoT}, $w$ enjoys the structure \eqref{structure:1Dsol} with $t\in[0,T]$. Recall the definition of $x_\f$ as in \eqref{estimate:x_infty} and $x_1$ as in \eqref{def:x_n}. Let $\underline{x},\overline{x}$ be defined as follows:
	$
	\underline{x}=x_0-MT\mbox{ and }\overline{x}=x_\f+MT.
	$ Now we want to show that $\abs{w}_{W^{s,p}([\underline{x},\overline{x}])}=\f$. In order to do so we define
	\begin{align}
	h_n&:={\e_n}{(T-t)}\mbox{ for }n\geq1,\label{def:h_n}\\
	\Omega_i[h,t]&:=\left\{z\in\R;\,x_{2i+1}+P_kt-h<z<x_{2i+1}+P_kt\right\}\mbox{ for }J_{k-1}<i\leq J_{k}.\nonumber
	\end{align}
	Note that if $h\in(h_{n+1},h_n]$ then for $J_{m-1}<i<J_{m}$ we have $
	x_{2i}+f^{\p}(a_k)t<x_{2i+1}+P_kt-h$ and $x_{2i+1}+P_kt+h<x_{2i+2}+P_kt$. Hence we obtain
	\begin{eqnarray}
	\|\De^{h}w(\cdot,t)\|_{L^p(\underline{x},\overline{x})}^p\geq\sum\limits_{k=1}^{n}\sum\limits_{i=J_{k-1}+1}^{J_k}\int\limits_{\Omega_i[h,t]}|a_k-a_0|^pdx &\geq& C_1\sum\limits_{k=1}^{n}hN_k|a_k-a_0|^p\nonumber\\
	&\geq& hn^{-pq_n+p_0q_n-1-\e}.\nonumber
	\end{eqnarray}
	Therefore, we have the following:
	\begin{align}
	|w(\cdot,t)|_{B^{s,p,\theta}([\underline{x},\overline{x}])}^\theta 
	&\geq\sum\limits_{n=n_0}^{\f}\int\limits_{h_{n+1}}^{h_n}\|\De_{\nu_0}^{h}w(\cdot,t)\|_{L^p(\underline{x},\overline{x})}^{\theta}\frac{1}{h^{1+s\theta}}dh\nonumber\\
	&\geq C_1\sum\limits_{n=n_0}^{\f}n^{-q_n\theta-\frac{\theta}{p}-\frac{\theta}{p}\e+\frac{p_0q_n\theta}{p}}h_n^{-s\theta}
	\int\limits_{h_{n+1}}^{h_n}h^{\frac{\theta}{p}-1}dh.\nonumber
	\end{align}
	By fundamental theorem of calculus we have
	\begin{align}
	|w(\cdot,t)|_{B^{s,p,\theta}([\underline{x},\overline{x}])}^\theta
	&\geq 
	 C_1\sum\limits_{n=1}^{\f}n^{-q_n\theta-\frac{\theta}{p}-\frac{\theta}{p}\e+\frac{p_0q_n\theta}{p}}h_n^{-s\theta+\frac{\theta}{p}}\left[1-\left(\frac{h_{n+1}}{h_n}\right)^{\frac{\theta}{p}}\right].\label{estimate:wn}
	\end{align}
	Recall definition of $\e_n$ as in \eqref{theorem2:asmp1}. From \eqref{def:h_n} we have $h_n=(T-t)\abs{a_k-a_0}^{p_0}$. From \eqref{theorem2:asmp2} we infer $h_n=(T-t)n^{-pq_n}$. Hence we have 
	\begin{equation}\label{estimate:h1}
	\frac{h_{n+1}}{h_n}=(n+1)^{-p_0q_{n+1}}n^{p_0q_n}=\left(\frac{n}{n+1}\right)^{p_0q_n}(n+1)^{p_0(q_n-q_{n+1})}.
	\end{equation}
	Recall from \eqref{WLOG:akqk} that $q_n\geq0$ and $q_{n+1}\geq q_n+1$, equivalently, $q_n-q_{n+1}\leq-1$. From \eqref{estimate:h1} we get
	$
	\frac{h_{n+1}}{h_n}\leq (n+1)^{p_0(q_n-q_{n+1})}\leq (n+1)^{-p_0}\leq \frac{1}{2^{p_0}}$ for $n\geq1.
	$
	Hence, from \eqref{estimate:wn} we obtain
	\begin{align}
	|w(\cdot,t)|_{B^{s,p,\theta}([\underline{x},\overline{x}])}^\theta
	&\geq  C_12^{-p}\sum\limits_{n=1}^{\f}n^{-q_n\theta-\frac{\theta}{p}-\frac{\theta}{p}\e+\frac{p_0q_n\theta}{p}}h_n^{-s\theta+\frac{\theta}{p}}.\label{estimate:wn2}
	\end{align}
	Using $h_n=(T-t)n^{-p_0q_n}$ in \eqref{estimate:wn2}, we infer
	\begin{align}
	|w(\cdot,t)|_{B^{s,p,\theta}([\underline{x},\overline{x}])}^\theta
	&\geq  C_12^{-1}\sum\limits_{n=1}^{\f}n^{-q_n\theta-\frac{\theta}{p}-\frac{\theta}{p}\e+\frac{p_0q_n\theta}{p}}n^{sp_0\theta q_n-\frac{\theta p_0q_n}{p}}.\label{estimate:wn3}
	\end{align}
	To simplify the exponent in \eqref{estimate:wn3} we observe the following
	\begin{equation*}
	-q_n\theta-\frac{\theta}{p}-\frac{\theta}{p}\e+\frac{p_0q_n\theta}{p}+sp_0\theta q_n-\frac{\theta p_0q_n}{p}=q_n\theta(sp_0-1)-\theta/p-\e\theta/p.
	\end{equation*}
	By choice of $p_0$ we have $sp_0>1$. Note that $q_n\rr\f$ as $n\rr\f$. Then there exists an $n_0$ such that
	$
	q_n\geq (sp_0-1)^{-1}\left[\theta/p+\e\theta/p-1\right].
	$ Hence we have $
	-q_n\theta-\frac{\theta}{p}-\frac{\theta}{p}\e+\frac{p_0q_n\theta}{p}+sp_0\theta q_n-\frac{\theta p_0q_n}{p}\geq -1$ for $n\geq n_0$. From \eqref{estimate:wn3} we infer
	\begin{align}\label{Wsp_blowup:w}
	|w(\cdot,t)|_{B^{s,p,\theta}([\underline{x},\overline{x}])}^\theta
	&\geq  C_12^{-1}\sum\limits_{n=n_0}^{\f}n^{-1}=\f.
	\end{align}
	As we observe in proof of Proposition \ref{Proposition:1d_uptoT} that $w_0$ is not a compactly supported data when $a_0\neq0$. We wish to find a data $v_0^T$ such that it has compact support and the corresponding entropy solution $v^T$ does not belong to $B^{s,p,\theta}$. Here we take the similar strategy as we have done in the proof of Proposition \ref{Proposition:1d_uptoT}. Let $M$ be defined as in \eqref{def:Wsp:M}.
	We consider $R_2=3MT$ and $I=[x_1-R_2,x_\f+R_2]$. We define $v_0^T$ as $v_0^T=w_0\chi_I$. By a similar analysis as we have done in proof of Proposition \ref{Proposition:1d_uptoT} we can show that
	\begin{equation}\label{nonlinear(Wsp):w=v}
	w(x,t)=v^T(x,t) \mbox{ for all }(x,t)\in [x_1-MT,x_\f+MT]\times[0,T]. 
	\end{equation}
	From \eqref{def:Wsp:M} we have $M\geq \abs{f^{\p}(a_0)}$. Therefore, we have
	\begin{equation}\label{nonlinear(Wsp):cal3}
	[\underline{x},\overline{x}]=[x_1-f^{\p}(a_0)T,x_\f+f^{\p}(a_0)T]\subset[x_1-MT,x_\f+MT].
	\end{equation}
	By \eqref{Wsp_blowup:w}, \eqref{nonlinear(Wsp):w=v} and \eqref{nonlinear(Wsp):cal3} we have $|v^T(\cdot,t)|_{B^{s,p,\theta}([\underline{x},\overline{x}])}=\f$ for all $ t\in[0,T]$. By a similar argument as in proof of Proposition \ref{Proposition:1d_uptoT} we find a $\la>0$ independent of $T$ such that $supp(v_0^T)\subset[-\la T,\la(T+1)]$ and $|v^T(\cdot,t)|_{B^{s,p,\theta}([-\la T,\la(T+1)])}=\f\mbox{ for all }t\in[0,T]$.	This concludes Proposition \ref{Proposition:1dW_uptoT}.
\end{proof}
\subsection{Proof of Theorem \ref{theorem2}}
Now we are ready to show Theorem \ref{theorem2}. Proof is similar to Theorem \ref{theorem1}. Here we use regularity property \eqref{mild_regularity} of flux and Proposition \ref{Proposition:1dW_uptoT} from section \ref{sec:Wsp1D}.
\begin{proof}[Proof of Theorem \ref{theorem2}]
	By our assumption $F\in C^2(I,\R^d)$ and satisfies \eqref{alpha_sup} for $\al_{sup}(F,I)<1$ along with mild regularity assumption \eqref{mild_regularity}. Fix an $s$ such that $\al_{sup}(F,I)<s<1$. Since $F$ satisfies \eqref{mild_regularity}, there exists a $\xi_0\in\mathcal{S}^{d-1}$ and $k_0\in\mathscr{I}_{\xi_0}$ such that $	\al_{sup}(\xi_0\cdot F,I_{k_0}^{\xi_0})<s$ where $I^{\xi_0}_{k_0}$ is defined in \eqref{intervals}. We define $f\in C^2(\re,\re)$ as $f=\xi_0\cdot F$. Observe that $f$ satisfies all the assumptions of Proposition \ref{Proposition:1dW_uptoT}. By a similar method as in proof of Theorem \ref{theorem1} we can construct an entropy solution $U_{\xi_0}$ such that $U_{\xi_0}(x,t)=v^T(x\cdot\xi_0,t)$ for $t\in[0,T]$. Note that
	$
	\|\De_i^hU_\xi(\cdot,t)\|_{L^p(B(\textbf{0},\la(T+1)))}=c_0	\|\De^hv^T(\cdot,t)\|_{L^p(B(\textbf{0},\la(T+1)))}.
	$  Hence we show that $\abs{U_{\xi_0}(\cdot,t)}_{B^{s,p,\theta}(B(\textbf{0},\la(T+1)))}=\f$ for $t\in[0,T]$. By a similar method as in Theorem \ref{theorem1} we construct an entropy solution $U$ such that $U(\cdot,t)\notin B^{s,p,\theta}_{loc}(\R^d)$ for all $t>0$. This finishes the proof.
\end{proof}

\subsection{Proof of Proposition \ref{proposition2}}

Now we can prove Proposition \ref{proposition2} in a similar manner as we did for Theorem \ref{theorem2}.
\begin{proof}[Proof of Proposition \ref{proposition2}]
	Proof is similar to Theorem \ref{theorem2}. In proof of Theorem \ref{theorem2}, we use condition \ref{mild_regularity} satisfied by flux $F$ to get direction $\xi$ and sequence $\{a_k\}$. Here flux need not satisfy \eqref{mild_regularity} but existence of direction $\xi$ and $\{a_k\}$ are guaranteed by Lemma \ref{lemma:Holder}.
\end{proof}


\subsection{Proof of Remark \ref{remark:BV_for_alpha=1}}\label{sec:Proof_of_Remark}
Combining Lemma \ref{lemma:logflux} and Theorem \ref{theorem:1d} we conclude Remark \ref{remark:BV_for_alpha=1}. 
\subsection{Proof of Proposition \ref{proposition1}}
Proposition \ref{proposition1} can be proved in a similar way as Theorem \ref{theorem2}.

\appendix\section{}

Now we prove Lemma \ref{lemma:planar_multiD}. It is an application of Proposition \ref{Prop:Kruzkov}.
\begin{proof}[Proof of Lemma \ref{lemma:planar_multiD}:] Fix a $\bar{\varphi}\in C^{\f}_{c}(\re^d\times[0,\f))$. Let $v,v_0$ be defined as in \eqref{lemma:kruzkov_v}. 
	Then we have 
	\begin{align*}
	& \int\limits_{0}^{\f}\int\limits_{\re^d}\left[v(x,t)\frac{\pa \bar{\varphi}}{\pa t}+\sum\limits_{i=1}^{d}f_i(v)\frac{\pa\bar{\varphi}}{\pa x_{i}}\right]dx dt+\int\limits_{\re^d}v_0(x)\bar{\varphi}(x,0)dx\\
	&=\int\limits_{0}^{\f}\int\limits_{\re^d}\left[w(x_1,t)\frac{\pa \bar{\varphi}}{\pa t}+\sum\limits_{i=1}^{d}f_i(w(x_1,t))\frac{\pa\bar{\varphi}}{\pa x_{i}}\right]dx dt+\int\limits_{\re^d}w_0(x_1)\bar{\varphi}(x,0)dx.
	\end{align*}
	As $f_i(w(x_1,t))$ is independent of $x_i$ for $i\neq1$, we have
	$
	\int\limits_0^{\f}\int\limits_{\re^d}f_i(w(x_1,t))\frac{\pa\bar{\varphi}}{\pa x_{i}}dx dt=0\mbox{ for }i\neq1.
	$
	Since $w$ is a weak solution to Cauchy problem \eqref{w1}, by using Fubini's Theorem, we get 
	\begin{align*}
	& \int\limits_{0}^{\f}\int\limits_{\re^d}\left[v(x,t)\frac{\pa \bar{\varphi}}{\pa t}+\sum\limits_{i=1}^{d}f_i(v)\frac{\pa\bar{\varphi}}{\pa x_{i}}\right]dx dt+\int\limits_{\re^d}v_0(x)\bar{\varphi}(x,0)dx\\
	&=\int\limits_{\re^{d-1}}\left[\int\limits_{0}^{\f}\int\limits_{\re}\left(w\frac{\pa \bar{\varphi}}{\pa t}+f_1(w)\frac{\pa\bar{\varphi}}{\pa x_{1}}\right)dx_1 dt+\int\limits_{\re}w_0(x_1)\bar{\varphi}(x,0)dx_1 \right]dx_2\cdots dx_d=0.
	\end{align*}
	Since $w$ is a weak solution to Cauchy problem \eqref{w1}, we obtain
	\begin{equation*}
	\int\limits_{0}^{\f}\int\limits_{\re^d}\left[v(x,t)\frac{\pa \bar{\varphi}}{\pa t}+\sum\limits_{i=1}^{d}f_i(v)\frac{\pa\bar{\varphi}}{\pa x_{i}}\right]dx dt+\int\limits_{\re^d}v_0(x)\bar{\varphi}(x,0)dx=0.
	\end{equation*}
	By a similar argument, we can show that $v$ satisfies Kruzkov entropy inequality. 
	This proves Lemma \ref{lemma:planar_multiD}.
\end{proof}
Now we prove Lemma \ref{lemma:linear_transform}. It follows from Proposition \ref{Prop:Kruzkov} after a suitable change of variable.
\begin{proof}[Proof of Lemma \ref{lemma:linear_transform}:]
	Fix a ${\varphi}\in C^{\f}_{c}(\re^d\times[0,\f))$. Since $\bar{v}$ is entropy solution to (\ref{eqn:conlaw}) with initial data $\bar{v}_0$ and flux $\bar{f}$, then we put $\bar{\varphi}(x)=\varphi(L^{-1}(x-\textbf{c}))\in C^{\f}_{c}(\re^d\times[0,\f))$ in weak formulation of $v$ and use $y=L^{-1}(x-\textbf{c})$ to obtain  
	\begin{align}
	&|\mbox{det}{L}|\int\limits_{0}^{\f}\int\limits_{\re^d}\bar{v}(L(y)+\textbf{c},t)\frac{\pa \bar{\varphi}}{\pa t}(L(y)+\textbf{c},t)+\sum\limits_{i=1}^{d}\bar{f}_i(\bar{v}(L(y)+\textbf{c},t))\frac{\pa\bar{\varphi}}{\pa x_{i}}(L(y)+\textbf{c},t)dy dt\nonumber\\
	&+|\mbox{det}{L}|\int\limits_{\re^d}\bar{v}_0(L(y)+\textbf{c})\bar{\varphi}(L(y)+\textbf{c},0)dy=0.\label{lemma:kruzkov2:cal1}
	\end{align}
	By definition of $\bar{\phi}$ and $\bar{F}$ we have
	$
	\frac{\pa\bar{\varphi}}{\pa x_{i}}(L(y)+\textbf{c},t)=\sum\limits_{j=1}^{d}(L^{-1})_{ij}\frac{\pa{{\varphi}}}{\pa x_{j}}(y,t)$ and $\sum\limits_{i=1}^{d}(L^{-1})_{ij} \linebreak \bar{f}_i(\bar{v}(L(y)+\textbf{c},t))=f_j(\bar{v}(L(y)+\textbf{c},t)).
     $
	Hence, from \eqref{lemma:kruzkov2:cal1} and definition of $\bar{v}_0$, we infer
	\begin{align*}
	&\int\limits_{0}^{\f}\int\limits_{\re^d}\bar{v}(L(y)+\textbf{c},t)\frac{\pa \varphi}{\pa t}(y,t)+\sum\limits_{j=1}^{d}f_j(\bar{v}(L(y)+\textbf{c},t))\frac{\pa{{\varphi}}}{\pa x_{j}}(y,t)dy dt+\int\limits_{\re^d}v_0(y){\varphi}(y,0)dy\\
	&=0.
	\end{align*}
	Since $\bar{v}$ satisfies the Kru\v{z}kov entropy inequality it is immediate that $\bar{v}(L(x)+\textbf{c},t)$ also satisfies the Kru\v{z}kov entropy inequality. Hence,  $\bar{v}(L(x)+\textbf{c},t)$ is an entropy solution to (\ref{eqn:conlaw}) with the initial data $v_0$. Hence from Proposition \ref{Prop:Kruzkov} we conclude that $v(x,t)=\bar{v}(L(x)+\textbf{c},t)$. This completes proof of Lemma \ref{lemma:linear_transform}.
\end{proof}

\noindent\textbf{Acknowledgement.} Authors acknowledge the support of the Department of Atomic Energy, Government of India, under project no. 12-R\&D-TFR-5.01-0520. The first author would like to thank Inspire faculty-research grant DST/INSPIRE/04/2016/00-0237. Authors thank IFCAM project ``Conservation laws: $BV^s$, interface and control".  
Authors are very thankful to  Pierre-Emmanuel Jabin for his valuable comments and suggestions. 

\end{document}